\documentclass{amsart}
 
\usepackage{xy}
\usepackage{amssymb}
\usepackage{amsmath}
\usepackage{amsthm} 
\usepackage[cmtip,all]{xypic}

\theoremstyle{plain}
\newtheorem{Lem}{Lemma}[section]

\newtheorem{Cor}[Lem]{Corollary}
\newtheorem{Thm}[Lem]{Theorem}

{\theoremstyle{definition} 
\newtheorem{Ex}[Lem]{Example}
\newtheorem{Rk}[Lem]{Remark}
\newtheorem{Def}[Lem]{Definition}}

\newcommand{\zig}{\addtocounter{Lem}{1}\tag{\theLem}}
\newcommand{\lhat}{\widehat{L}}
\def\:{\colon}
\DeclareMathOperator*{\colim}{colim}
\DeclareMathOperator*{\holim}{holim}

\DeclareMathOperator*{\cdp}{cd_\mathit{p}}
\DeclareMathOperator*{\cdq}{cd_\mathit{q}}

\DeclareMathOperator*{\Zq}{\mathbb{Z}_\mathit{q}}

\DeclareMathOperator*{\GG}{\mathbb{G}}
\DeclareMathOperator*{\Zp}{\mathbb{Z}/\mathit{p}}
\DeclareMathOperator*{\gpring}{\mathbb{Z}/\mathit{p}[\mathbb{Z}/\mathit{q}^\mathit{n}]}
\DeclareMathOperator*{\Hn}{\Sigma^\mathit{n}\mathit{H}(\mathbb{Z}/\mathit{p}[\mathbb{Z}/\mathit{q}^\mathit{n}])}

\begin{document}


\title{Iterated homotopy fixed points for the Lubin-Tate spectrum}          

\author[Daniel G. Davis]{Daniel G. Davis$\sp 1$}

\email{dgdavis@louisiana.edu}
\address{Department of Mathematics, University of Louisiana at Lafayette, Lafayette, LA 70504 
\vspace{.1in}}

\email{wieland@math.brown.edu} 
\address{Department of Mathematics, Brown University, Providence, RI 02912}

\footnotetext[1]{The author was partially supported by an NSF VIGRE grant at Purdue University. Part of this 
paper was written during a visit to the Institut Mittag-Leffler (Djursholm, 
Sweden) and a year at Wesleyan University.}

\maketitle
\let\languagename\relax

\vspace{-0.5cm}
\begin{center}
\textsc{with an appendix by daniel g. davis$\sp 2$ and ben wieland$\sp 3$}
\end{center}
\footnotetext[2]{The author was supported by a grant from the Louisiana Board 
of Regents Support Fund.}
\footnotetext[3]{The author was partially supported by an NSF Postdoctoral Fellowship.}

\begin{abstract} 
When $G$ is a profinite group and $H$ and $K$ are closed subgroups, with 
$H$ normal in $K$, it is not known, in general, how to form the iterated homotopy 
fixed point spectrum $(Z^{hH})^{hK/H},$ where $Z$ is a 
continuous $G$-spectrum and all group actions are to be continuous. 
However, we show that, if $G=G_n$, the 
extended Morava stabilizer group, and $Z=\lhat(E_n \wedge X),$ 
where $\lhat$ is Bousfield localization with respect to Morava $K$-theory, 
$E_n$ is the Lubin-Tate spectrum, and $X$ is any spectrum with trivial 
$G_n$-action, then the iterated homotopy fixed point spectrum can always 
be constructed. Also, we 
show that $(E_n^{hH})^{hK/H}$ is just $E_n^{hK}$, extending a result of 
Devinatz and Hopkins.
\end{abstract}

\section{Introduction}
\par
Let $G$ be a profinite group and let $X$ be a continuous $G$-spectrum. Thus, 
$X$ is the homotopy limit of a tower 
\[X_0 \leftarrow X_1 \leftarrow X_2 \leftarrow \cdots \leftarrow X_i \leftarrow \cdots \]
of discrete $G$-spectra that are 
fibrant as spectra (see Section \ref{summary} for a quick review of this 
notion that is developed in detail in \cite{cts}). 
Let $H$ and $K$ be closed subgroups of $G$, with 
$H$ normal in $K$. Note that $H$ is closed in $K$, and $K$, $H$, and 
the quotient $K/H$ are all profinite groups. 
\par
It is easy to see that the fixed point spectrum $X^H$ is a $K/H$-spectrum 
and \[(X^{H})^{K/H} = X^K.\] Now, there is a model category $\mathcal{T}$ 
of towers of discrete 
$G$-spectra and homotopy fixed points 
are the total right derived functor of 
\[\lim_i (-)^G \: \mathcal{T} \rightarrow \mathrm{spectra}\] (see \cite[Section 4, Remark 8.4]{cts} for 
more detail). Thus, as noted in \cite[Introduction]{LHS}, 
since homotopy fixed points are defined in terms of fixed points, 
it is reasonable to wonder (a) if the $H$-homotopy 
fixed point spectrum 
$X^{hH}$ is a continuous $K/H$-spectrum; and (b) assuming that (a) 
holds, if there is an equivalence 
\begin{equation}\zig\label{question}
(X^{hH})^{hK/H} \simeq X^{hK},
\end{equation} where $(X^{hH})^{hK/H}$ is the {\em iterated homotopy fixed point 
spectrum}. 
\par
It is not hard to see that, if $H$ is open in $K$, then 
the isomorphism of (\ref{question}) always holds; see 
Theorem \ref{finite} for the details. However, if $H$ is not open in $K$, 
then the situation is much more complicated. In this case, 
it is not known, in general, how to construct $X^{hH}$ as a $K/H$-spectrum, 
and, even when we can construct $X^{hH}$ as a $K/H$-spectrum, it is not known, 
in general, how to view $X^{hH}$ as a continuous $K/H$-spectrum. Thus, when 
$H$ is not open in $K$, additional hypotheses are needed just to get past 
step (a)  - see 
Sections \ref{details} and \ref{hyperfibrant} for a discussion of this point.
\par
Now we consider the above issues in an example that is of interest in 
chromatic stable homotopy theory. 
Let $n \geq 1$, let $p$ be a fixed prime, and let $K(n)$ be the $n$th 
Morava $K$-theory spectrum (so that $K(n)_\ast = \mathbb{F}_p[v_n^{\pm 1}]$, 
where $|v_n| = 2(p^n-1)$). Then, let $E_n$ be the Lubin-Tate spectrum: 
$E_n$ is the $K(n)$-local Landweber exact 
spectrum whose coefficients are given by 
$E_{n \ast} = 
W(\mathbb{F}_{p^{n}})[[u_1, ..., u_{n-1}]][u^{\pm 1}],$ where 
$W(\mathbb{F}_{p^{n}})$ is the ring of Witt vectors of the field 
$\mathbb{F}_{p^{n}}$, each $u_i$ has degree zero, and the degree of $u$ is 
$-2$. 
\par
Let $S_n$ be the $n$th Morava stabilizer group (the automorphism 
group of the Honda formal group law $\Gamma_n$ of height $n$ over 
$\mathbb{F}_{p^n}$) and let \[G_n = 
S_n \rtimes \mathrm{Gal}(\mathbb{F}_{p^{n}}/\mathbb{F}_p)\] be 
the extended Morava 
stabilizer group, a profinite group. Then, given any closed 
subgroup $K$ of $G_n$, 
let $E_n^{dhK}$ be the construction, due to Devinatz and Hopkins, that 
is denoted $E_n^{hK}$ in \cite{DH}. This work (\cite{DH}) shows that 
the spectra $E_n^{dhK}$ behave like homotopy fixed point spectra, though 
they are not constructed with respect to a continuous $K$-action on $E_n$.
\par
As in \cite{cts}, let $E_n^{hK}$ be the $K$-homotopy fixed point spectrum of $E_n$, 
formed with respect to the continuous $K$-action on $E_n$. In \cite{thesis} 
(see \cite[Theorem 8.2.1]{joint} 
for a proof), the author showed that 
\begin{equation}\zig\label{thesis}
E_n^{hK} \simeq E_n^{dhK},
\end{equation} 
for all $K$ closed in $G_n$. 
Now, as explained in \cite[pg. ~8]{LHS}, the profinite group $K/H$ 
acts on $E_n^{dhH}$. Suppose that $H$ is open in $K$, so 
that $K/H$ is a finite group. Then, 
Devinatz and Hopkins (see \cite[Theorem 4 and Section 7]{DH}) showed that the 
canonical map 
\[E_n^{dhK} \rightarrow \holim_{K/H} E_n^{dhH} = 
(E_n^{dhH})^{hK/H}\] 
is a weak equivalence. Thus, by applying (\ref{thesis}),
we obtain that 
$(E_n^{hH})^{hK/H} \simeq E_n^{hK}$. Alternatively, as mentioned earlier, 
since $E_n$ is a continuous $K$-spectrum 
and $H$ is open in $K$, this same conclusion follows immediately from 
Theorem \ref{finite}. Therefore, it is 
clear that, for $E_n$ and the finite 
quotients $K/H$, the iterated homotopy fixed point spectrum behaves 
as desired - that is, the equivalence of (\ref{question}) is valid.
\par
However, when $K/H$ is not finite, the papers \cite{DH, LHS} do not 
construct $(E_n^{dhH})^{hK/H}$, and hence, they are not 
able to consider the question of whether this spectrum is just 
$E_n^{dhK}$. Thus, in this more complicated situation, we 
show how to construct $(E_n^{dhH})^{hK/H}$. More precisely, we prove the 
following result, where $\lhat$ denotes the Bousfield localization functor 
$L_{K(n)}$.
\begin{Thm}\label{main}
Let $H$ and $K$ be any closed subgroups of $G_n$, such that $H$ is a normal 
subgroup of $K$, so that $E_n^{dhH}$ carries the natural $K/H$-action. Also, 
let $X$ be any spectrum $\mathrm{(}with$ no $K/H$-action$\mathrm{)}$. Then the spectrum 
$\lhat(E_n^{dhH} \wedge X)$ is a continuous $K/H$-spectrum, with $K/H$ acting 
diagonally on $\lhat(E_n^{dhH} \wedge X)$ $\mathrm{(}$$X$ has the trivial 
$K/H$-action$\mathrm{)}$.
\end{Thm}
\par
The theory of \cite{cts} (reviewed in Section 2), applied to Theorem 
\ref{main}, automatically yields 
\[(\lhat(E_n^{dhH} \wedge X))^{hK/H}\] as the homotopy fixed points with 
respect to the continuous $K/H$-action. Running through the argument for 
Theorem \ref{main} again, by omitting the $(- \wedge X)$ everywhere, 
gives the desired object 
$(E_n^{dhH})^{hK/H}.$ By applying (\ref{thesis}), we immediately 
find that the iterated homotopy fixed point spectrum 
$(E_n^{hH})^{hK/H}$ is always defined.
\par
In Section \ref{dsssection}, we use the descent spectral sequence for 
$(\lhat(E_n^{dhH} \wedge E_n))^{hK/H}$, where $K/H$ is acting trivially on (the second) $E_n$, to prove the following result, 
which says that the iterated homotopy fixed point spectrum for $E_n$ behaves 
as in (\ref{question}), \emph{if} one smashes with $E_n$ before taking the $K/H$-homotopy 
fixed points. 
\begin{Thm}
If $H$ and $K$ are as in Theorem \ref{main}, then 
\[(\lhat(E_n^{dhH} \wedge E_n))^{hK/H} \simeq
\lhat(E_n^{dhK} \wedge E_n).\]
\end{Thm}
\par
In Section \ref{morava}, we use (\ref{thesis}) to prove that the iterated 
homotopy fixed point spectrum for $E_n$ behaves as desired.
\begin{Thm}\label{lastmain}
Let $H$ and $K$ be as above. Then there is an equivalence 
\[(E_n^{hH})^{hK/H} \simeq E_n^{hK}.\]
\end{Thm}
\par
These iterated homotopy fixed point spectra, $(E_n^{hH})^{hK/H}$, play a useful role in 
chromatic homotopy theory. A major conjecture 
in this field is that $\pi_\ast(\lhat(S^0))$ is a module of finite type over the $p$-adic integers 
$\mathbb{Z}_p$. An important way to tackle this conjecture is due to Devinatz 
(see \cite{homotopydev, finite}), and a key part of his program is to find a 
closed subgroup $H_0$ of $G_n$ and a finite spectrum $Z$ that is not rationally acyclic, such that 
$\pi_\ast(E_n^{hH_0} \wedge Z)$ is of finite type and $H_0$ is a part of a chain 
\[H_0 \vartriangleleft H_1 \vartriangleleft \cdots \vartriangleleft H_t = G_n\] of closed 
subgroups. 
\par
For suppose that, for some $i < t$, $\pi_\ast(E_n^{hH_i} \wedge Z)$ is of finite type. 
By \cite[(0.1)]{LHS} and Theorem \ref{sscomparison}, 
there is a strongly convergent descent spectral sequence
\[H^s_c(H_{i+1}/H_i; \pi_t(E_n^{hH_i} \wedge Z)) \Rightarrow \pi_{t-s}((E_n^{hH_i})^{hH_{i+1}/H_i} 
\wedge Z),\] where the $E_2$-term is continuous cohomology for profinite 
$\mathbb{Z}_p[[H_{i+1}/H_i]]$-modules (defined just above (\ref{dss3})). As explained in \cite[pg.~133]{homotopydev}, various 
properties of these cohomology groups and the spectral sequence imply that 
\[\pi_\ast((E_n^{hH_i})^{hH_{i+1}/H_i} 
\wedge Z) \cong \pi_\ast(E_n^{hH_{i+1}} \wedge Z)\] has finite type. Thus, induction shows 
that $\pi_\ast(E_n^{hG_n} \wedge Z) \cong \pi_\ast(\lhat(Z))$ is of finite type, since 
$E_n^{hG_n} \simeq \lhat(S^0)$ (by \cite[Theorem 1(iii)]{DH}, 
\cite[Corollary 8.1.3]{joint}). In particular, since 
$Z$ is not rationally acyclic, it follows that $\pi_\ast(\lhat(S^0))$ would then have 
finite type (see \cite[pg. 2]{finite}).
\par
Other examples of the utility of the $K/H$-action on $E_n^{hH},$ with $K/H$ an infinite 
profinite group, are in 
\cite[Section 8]{DH} and \cite[Section 2]{GHM}: in the first case, the action helps (see 
\cite[Proposition 8.1]{DH}) to construct 
an interesting element in $\pi_{-1}(\lhat(S^0))$, for all $n$ and $p$ (\cite[Theorem 6]{DH}), and, 
in the second case, the action plays a role (see \cite[Theorem 10]{GHM}) 
in computing $\pi_\ast(L_2V(1))$, at the prime $3$.
\par
To study the iterated homotopy fixed points of $E_n$, we consider the notion of 
a hyperfibrant discrete $G$-spectrum, and we show that, for such a spectrum, 
the iterated homotopy fixed point spectrum is 
always defined in a natural way (see Definitions \ref{hyper} and 
\ref{construction}). Also, Lemma \ref{iterated} shows that if a discrete 
$G$-spectrum $X$ is a hyperfibrant discrete $K$-spectrum, then, for any $H$ closed in $G$ and 
normal in $K$, (\ref{question}) is valid. Additionally, we show 
that, for a totally hyperfibrant discrete $G$-spectrum, (\ref{question}) holds for all $H$ and 
$K$ closed in $G$, with $H$ normal in $K$ (see Definition \ref{totallyhyper}).  
\par
Given distinct primes $p$ and $q$, for $G = \mathbb{Z}/p \times \mathbb{Z}_q$, 
Ben Wieland has found an example of 
a discrete $G$-spectrum $X$ that is not a hyperfibrant 
discrete $G$-spectrum. In the Appendix, Wieland and the author provide the details for 
this example. As explained in Section \ref{details} (after Theorem \ref{finite}), 
for discrete $G$-spectra that are not hyperfibrant, there are situations where it is not known 
how to define a continuous $K/H$-action on $X^{hH}$, so that it is also not known 
how to define the $K/H$-homotopy fixed points of $X^{hH}$, with respect to a 
continuous $K/H$-action.  
\par
At the end of Section \ref{morava}, we point out that the results in this paper apply to 
certain other spectra that are like $E_n$, in that they replace the role of $\mathbb{F}_{p^n}$ and 
the height $n$ Honda formal group law $\Gamma_n$ with certain other finite fields and 
height $n$ formal group laws, respectively.   
\par
In \cite{joint}, Mark Behrens and the author show that if a spectrum $E$ is a consistent profaithful 
$k$-local profinite $G$-Galois extension of $A$ of finite vcd, then the map $\psi(E)_K^G$ of 
(\ref{sigh}) is a $k_\ast$-equivalence, for all $K$ closed in $G$ (\cite[Theorem 7.1.1]{joint}), where 
$L_k(-) \simeq L_ML_T(-)$, with $M$ a finite spectrum and $T$ smashing, and $A$ is a 
$k$-local commutative $S$-algebra. Thus, by \cite[Corollary 7.1.3]{joint}, if $K$ is normal in $G$, $E^{hK}$ is $k$-locally a 
discrete $G/K$-spectrum, so that, in a $k$-local context, it is natural to define $(E^{hK})^{hG/K} 
= (\colim_{N \vartriangleleft_o G} E^{hNK})^{hG/K}$, and, then, $L_k((E^{hK})^{hG/K}) \simeq 
L_k(E^{hG}).$ 
\par
We point out that, independently of our work, in 
\cite[Sections 10.1, 10.2]{fausk}, Halvard 
Fausk considers the iterated homotopy fixed point pro-spectrum for pro-G-orthogonal spectra, and the associated descent spectral sequence, by making general use of Postnikov towers. However, our approach to iterated homotopy fixed points is different in technique from his: instead of Postnikov towers, we rely on the notion of hyperfibrancy, which, in the case of $E_n$, makes use of technical results about the Morava stabilizer group (see \cite[the proof of Theorem 4.3]{DH}).  
\par
In \cite[Lemma 10.5]{wilkerson}, it is shown that, if 
$G$ is any discrete group, with $K$ any normal subgroup, and, if $X$ is a $G$-space, 
then $X^{hK}$ is homotopy equivalent to $(X^{hK})^{hG/K}$. This result, described as verifying a ``transitivity property" of homotopy fixed points, is, as 
far as the author has been able to determine, the 
first reference to iterated homotopy fixed points in the literature. Given an $S$-module 
$F$ and a faithful $F$-local 
$G$-Galois extension $A \rightarrow B$ of commutative $S$-algebras, where $G$ is a stably 
dualizable group, \cite[Theorem 7.2.3]{Rognes} shows that, if $K$ is an 
allowable normal subgroup of $G$, then $(B^{hK})^{hG/K} \simeq B^{hG}.$ 
\par
We make a comment about notation in this paper. 
If a limit or colimit is indexed 
over $N$, as in $\lim_N G/N$, then the (co)limit is indexed over all 
open normal subgroups of $G$, unless stated otherwise.     
\vspace{.1in}
\par
\noindent
\textbf{Acknowledgements.} I am grateful to Mark Behrens for a series of profitable 
conversations about iterated homotopy fixed points. I thank Tilman Bauer, Ethan Devinatz, and Paul Goerss for helpful discussions about this work. Also, I am grateful to the referee for many helpful comments that improved and sharpened the writing of this paper. Additionally, 
I thank Halvard Fausk, Rick Jardine, Ben Wieland, and an earlier referee for their comments. 
\section{A summary of the theory of continuous $G$-spectra and 
their homotopy fixed points}\label{summary}
This section contains a quick review of material from \cite{cts} that is needed for our 
work in this paper. We begin with some terminology. 
\par
All of our spectra are 
Bousfield-Friedlander spectra of simplicial sets. Let $G$ be a profinite group. A {\em discrete 
$G$-set} is a $G$-set $S$ such that the action map 
$G \times S \rightarrow S$ is continuous, where $S$ is regarded 
as a discrete space.
\begin{Def}
A {\em discrete $G$-spectrum} $X$ is a $G$-spectrum such that each 
simplicial set $X_k$ is a simplicial object in the category of discrete 
$G$-sets. $\mathrm{Spt}_G$ is the category of discrete $G$-spectra, 
where the morphisms are 
$G$-equivariant maps of spectra.
\end{Def}
\begin{Thm}[{\cite[Theorem 3.6]{cts}}]
$\mathrm{Spt}_G$ is a model category, where $f$ in $\mathrm{Spt}_G$ is 
a weak equivalence $\mathrm{(}$cofibration$\mathrm{)}$ if and only 
if it is a weak equivalence $\mathrm{(}$cofibration$\mathrm{)}$ of spectra.
\end{Thm}
\par
It is helpful to note that, by \cite[Lemma 3.10]{cts}, if $X$ is fibrant in $\mathrm{Spt}_G$, then $X$ is fibrant as a spectrum. In this paper, we often take the smash product (in spectra) of a discrete $G$-spectrum $X$ with a spectrum $Y$ that has no $G$-action (and we never take the smash product in the case where $Y$ has a non-trivial $G$-action). Then (following \cite[Lemma 3.13]{cts}), since $X^N$, where $N$ is an open normal subgroup of $G$, is a $G/N$-spectrum, the isomorphisms
\[X \wedge Y \cong (\colim_N X^N) \wedge Y \cong \colim_N (X^N \wedge Y)\] show that 
$X \wedge Y$ is also a discrete $G$-spectrum.  
\begin{Def}
If $X \in \mathrm{Spt}_G,$ then let 
$X \rightarrow X_{f,G} \rightarrow \ast$ be a trivial cofibration 
followed by a fibration, all in $\mathrm{Spt}_G$. Then 
$X^{hG} = (X_{f,G})^G.$ 
\end{Def}
\par
When $G$ is a finite group, this definition is equivalent to the 
classical definition of $X^{hG}$ for a finite discrete group $G$ (see \cite[Section 5]{cts}).
\begin{Def}
A {\em tower of discrete $G$-spectra} $\{X_i\}$ is a diagram 
\[X_0 \leftarrow X_1 \leftarrow X_2 \leftarrow \cdots\] 
in $\mathrm{Spt}_G$, such that each $X_i$ is fibrant as a spectrum. 
\end{Def}
\begin{Def}\label{hg}
If $\{X_i\}$ is a tower of discrete $G$-spectra, then $\holim_i X_i$ 
is a {\em continuous $G$-spectrum}. Also, $\holim_i (X_i)^{hG}$ is the 
{\em $G$-homotopy fixed point spectrum} of the continuous $G$-spectrum 
$\holim_i X_i$. Given a tower $\{X_i\}$, we write $X = \holim_i X_i$ 
and $X^{hG} = \holim_i (X_i)^{hG}.$
\end{Def}
\par
We call the above construction ``homotopy fixed points'' because: 
(1) if $G$ is a finite discrete group, then the above definition 
agrees with the classical definition of $X^{hG}$ (see \cite[Lemma 8.2]{cts}); and (2) it is the 
total right derived functor of fixed points in the appropriate sense (see \cite[Remark 8.4]{cts}). 
\begin{Def}
We say that $G$ has {\em finite virtual cohomological dimension} 
({\em finite vcd}), and we write $\mathrm{vcd}(G) < \infty$, if there 
exists an open subgroup $U$ in $G$ such that, for some positive integer $m$, 
$H^s_c(U;M) = 0$, for all $s > m$, whenever $M$ is a discrete $U$-module.
\end{Def}
\par
It is useful to note that if $G$ is a compact $p$-adic analytic group, 
then $G$ has finite vcd (see the explanation 
just before Lemma 2.10 in \cite{cts}).
\par
Before stating the next result, we need some notation. When $X$ is a discrete set, let 
$\mathrm{Map}_c(G,X)$ denote the discrete $G$-set of continuous maps $G \rightarrow X$, where the $G$-action is given by $(g \cdot f)(g') = f(g'g),$ for $g, g' \in G$. 
Now consider the functor 
\[\Gamma_G \colon \mathrm{Spt}_G \rightarrow \mathrm{Spt}_G, \ \ \ 
X \mapsto \Gamma_G(X) = \mathrm{Map}_c(G,X),\] where the action of $G$ on 
$\mathrm{Map}_c(G,X)$ 
is induced on the level of sets 
by the $G$-action on $\mathrm{Map}_c(G,(X_k)_l),$ the set of $l$-simplices of the pointed 
simplicial set $\mathrm{Map}_c(G,X_k),$ where $k, l \geq 0.$ 
\par
Let $\mathrm{Spt}$ denote the category of spectra. Then 
the right adjoint to the forgetful functor $U \colon \mathrm{Spt}_G \rightarrow \mathrm{Spt}$ is the 
functor $\mathrm{Map}_c(G, - )$. Also, as explained 
in \cite[Definition 7.1]{cts}, the functor $\Gamma_G$ 
forms a triple and there is 
a cosimplicial discrete $G$-spectrum $\Gamma^\bullet_G X.$ Given a tower 
$\{M_i\}$ of discrete $G$-modules, let $H^s_\mathrm{cont}(G;\{M_i\})$ denote 
continuous cohomology in the sense of Jannsen (see \cite{Jannsen}).
\begin{Thm}[{\cite[Theorem 8.8]{cts}}]\label{dss}
If $G$ has finite vcd and $\{X_i\}$ is a tower of discrete 
$G$-spectra, then there is a conditionally convergent 
descent spectral sequence 
\[E_2^{s,t} \Rightarrow \pi_{t-s}(\holim_i (X_i)^{hG}) = 
\pi_{t-s}(X^{hG}),\] where 
$E_2^{s,t} = \pi^s\pi_t(\holim_i(\Gamma^\bullet_G(X_i)_{f,G})^G).$ 
If the tower of abelian groups $\{\pi_t(X_i)\}$ satisfies the 
Mittag-Leffler condition for each $t \in \mathbb{Z}$, then 
\[E_2^{s,t} \cong H^s_\mathrm{cont}(G;\{\pi_t(X_i)\}).\]
\end{Thm}  
\section{Iterated homotopy fixed points - a discussion}\label{details}
\par
Let $G$ be a profinite group, and let $X$ be a continuous $G$-spectrum, 
so that $X = \holim_i X_i$, where each 
$X_i$ is a discrete $G$-spectrum that is 
fibrant as a spectrum. (Every discrete $G$-spectrum $Z$ is a continuous 
$G$-spectrum, in the sense that there is a weak equivalence $Z \rightarrow 
\holim_i Z_{f,G}$ that is given by the composition
\[Z \overset{\simeq}{\longrightarrow} Z_{f,G} 
\overset{\cong}{\longrightarrow} \lim_i Z_{f,G}  
\overset{\simeq}{\longrightarrow} \holim_i Z_{f,G},\] 
which is $G$-equivariant, 
where $\holim_i Z_{f,G}$ is a continuous $G$-spectrum 
and the limit and holim are each applied to a constant diagram.)
In this section, we consider more carefully the iterated homotopy fixed point spectrum 
$(X^{hH})^{hK/H}$ and its relationship with 
$X^{hK}$.  
\par
The author thanks Mark Behrens and Paul Goerss for help with the next 
lemma and its proof.
\begin{Lem}\label{forget}
Let $G$ be a profinite group and let $H$ be an open subgroup of 
$G$. If $X$ is a fibrant discrete $G$-spectrum, then $X$ is a fibrant 
discrete $H$-spectrum.
\end{Lem}
\begin{Rk}
If $H$ is normal in $G$, then Lemma \ref{forget} is a consequence of 
\cite[Remark 6.26]{Jardine}, using the fact that the presheaf of spectra 
$\mathrm{Hom}_G(-,X)$ is globally fibrant in the model category of presheaves 
of spectra on the site $G \negthinspace - \negthinspace \mathbf{Sets}_{df}$ 
(for an explanation of this fact, see \cite[Section 3]{cts}).
\end{Rk}
\begin{proof}[Proof of Lemma \ref{forget}.]
Note that 
the forgetful functor $U \: \mathrm{Spt}_G \rightarrow \mathrm{Spt}_H$ has 
a left adjoint 
$\mathrm{Ind}_G^H \: \mathrm{Spt}_H \rightarrow \mathrm{Spt}_G$, where, 
as a spectrum, $\mathrm{Ind}_G^H(Z)$ can be identified with the finite wedge 
$\bigvee_{G/H} Z$ of copies of $Z$ indexed by $G/H$.
\par
Let $X$ be a fibrant discrete $G$-spectrum. To show that $X$ is fibrant in 
$\mathrm{Spt}_H$, we will show that the forgetful functor $U$ preserves all 
fibrations. To do this, it suffices to show that $\mathrm{Ind}_G^H$ 
preserves all weak equivalences and cofibrations. Thus, we only have to 
show that if $f \: Y \rightarrow Z$ is a map of discrete $H$-spectra that is 
a weak equivalence (cofibration) of spectra, then 
$\mathrm{Ind}_G^H(f)$ is also a weak equivalence (cofibration) of spectra.
\par
If $f$ is a weak equivalence, then, as a map of spectra, 
$\pi_\ast(\mathrm{Ind}_G^H(f))$ is equivalent to $\oplus_{G/H} \pi_\ast(f),$ 
which is an isomorphism, so that $\mathrm{Ind}_G^H(f)$ is a weak equivalence. 
\par
Now let $f$ be a cofibration of spectra. Given $n \geq 0$, let, for example, $Y_n$ be the $n$th pointed simplicial set of $Y$ (a Bousfield-Friedlander spectrum of simplicial sets). 
Then $f_0 \: Y_0 \rightarrow Z_0$ 
and each induced map $(S^1 \wedge Z_n) \cup_{(S^1 \wedge Y_n)} Y_{n+1} 
\rightarrow Z_{n+1}$ are cofibrations of simplicial sets. It is easy to see 
that \[(\mathrm{Ind}_G^H(f))_0 \: \bigvee_{G/H} Y_0 \rightarrow 
\bigvee_{G/H} Z_0\] is a cofibration of simplicial sets. Also, the 
composition 
\[(S^1 \wedge \negthinspace \bigvee_{G/H} \negthinspace \negthinspace 
Z_n) \cup_{(S^1 \wedge 
\bigvee_{G/H} Y_n)} 
(\negthinspace \bigvee_{G/H} \negthinspace Y_{n+1}) \cong \negthinspace 
\bigvee_{G/H} \negthinspace ((S^1 \wedge Z_n) 
\cup_{(S^1 \wedge Y_n)} Y_{n+1}) \rightarrow \negthinspace \bigvee_{G/H} 
\negthinspace \negthinspace Z_{n+1}\] is a 
cofibration of simplicial sets. Thus, $\mathrm{Ind}_G^H(f)$ is a cofibration 
of spectra.
\end{proof}
\par
Let $G$ and $H$ be as in Lemma \ref{forget}. If $X$ is a discrete 
$G$-spectrum, then Lemma \ref{forget} implies that 
\begin{equation*}\label{corollary}\zig
X^{hH} = (X_{f,G})^H,
\end{equation*} since $X \rightarrow X_{f,G}$ is 
a trivial cofibration in $\mathrm{Spt}_H$.
\par
For the rest of this section, let $H$ and $K$ be closed subgroups of $G$, with $H$ normal in $K$. The following result shows that, if $H$ is open in $K$, then 
the equivalence of (\ref{question}) always holds.
\begin{Thm}\label{finite}
Let $G$ be a profinite group with closed subgroups $H$ and $K$, with $H$ 
normal in $K$. Given $\{X_i\}$, a tower of discrete $G$-spectra, let $X = \holim_i X_i$ 
be a continuous $G$-spectrum. If $H$ is open in $K$, then there 
is a weak equivalence \[X^{hK} \rightarrow (X^{hH})^{hK/H}.\]
\end{Thm}
\begin{proof}
Since each $X_i$ is a discrete $H$- and $K$-spectrum, $X$ is a continuous 
$H$- and $K$-spectrum. By (\ref{corollary}), 
\[X^{hH} = \holim_i (X_i)^{hH} = \holim_i ((X_i)_{f,K})^H.\] Since 
$\mathrm{Hom}_K(-,(X_i)_{f,K})$ is a fibrant presheaf of spectra on the 
site $K \negthinspace - \negthinspace \mathbf{Sets}_{df}$, the canonical map 
\[(X_i)^{hK} = \mathrm{Hom}_K(\ast,(X_i)_{f,K}) \rightarrow 
\holim_{K/H} \mathrm{Hom}_K(K/H,(X_i)_{f,K}) \cong 
\holim_{K/H} ((X_i)_{f,K})^H\] is a weak equivalence, by 
\cite[Proposition 6.39]{Jardine}. By (\ref{corollary}), there is the equality 
$\holim_{K/H} ((X_i)_{f,K})^H = \holim_{K/H} (X_i)^{hH},$ so that the previous 
sentence yields a weak equivalence 
\[(X_i)^{hK} \rightarrow \holim_{K/H} (X_i)^{hH}.\] Note that both the 
source and target of this weak equivalence are fibrant spectra, since, 
for any profinite group $L$, the 
functor $(-)^L \: \mathrm{Spt}_L \rightarrow \mathrm{Spt}$ preserves 
fibrant objects (see \cite[Corollary 3.9]{cts}).
\par
Thus, there is a weak equivalence 
\[X^{hK} = \holim_i (X_i)^{hK} \rightarrow \holim_i 
\holim_{K/H} (X_i)^{hH} \cong \holim_{K/H} X^{hH}.\] 
Since the identity map $X^{hH} \rightarrow X^{hH}$, where $X^{hH} = 
\holim_i ((X_i)_{f,K})^{H}$ is, of course, a weak equivalence and a 
$K/H$-equivariant map, with $X^{hH}$ fibrant, we can make the 
identification
$\holim_{K/H} X^{hH} = (X^{hH})^{hK/H},$ completing the proof. 
\end{proof} 
\par
It will be useful to define the following map. 
Given a discrete $G$-spectrum $X$, the fibrant replacement
$X_{f,G} \rightarrow (X_{f,G})_{f,K}$ and the identity 
\[\colim_{N} X^{hNK} = \colim_{N} (X_{f,G})^{NK},\] 
which applies (\ref{corollary}), induces the map
\begin{equation}\label{sigh}\zig
\psi(X)_K^G \: \colim_{N} X^{hNK} \cong (X_{f,G})^K 
\rightarrow ((X_{f,G})_{f,K})^K = X^{hK},
\end{equation} 
where the last equality is due to 
the fact that $X \rightarrow X_{f,G} \rightarrow (X_{f,G})_{f,K}$ is a 
trivial cofibration in $\mathrm{Spt}_K$, with $(X_{f,G})_{f,K}$ fibrant in 
$\mathrm{Spt}_K$. For convenience, we will sometimes write $\psi(X)_K^G$ as just $\psi$.
\par
Now we consider the case when (\ref{question}) is more difficult to understand; for the rest of 
this section, unless stated otherwise, we assume that $H$ is not open in $K$. Also, we assume 
that $H \neq \{e\}$, for this case is not problematic, since we have 
$(X^{h\{e\}})^{hK/\{e\}} = ((X_{f,G})^{\{e\}})^{hK} \simeq 
X^{hK}$. 
Additionally, we 
simplify the situation some by assuming that the continuous $G$-spectrum $X$ is 
a discrete $G$-spectrum. 
\par
First of all, in this case, consider the statement that 
$X_{f,K}$ is a fibrant object in $\mathrm{Spt}_H$. 
When $H$ is open in $K$, this statement is Lemma \ref{forget}, which we proved by 
showing that the forgetful functor 
$U \: \mathrm{Spt}_K \rightarrow \mathrm{Spt}_H$ preserves fibrations, and we did this by using the fact that $U$ is a right adjoint. However, in general, when $H$ is not open in $K$, the forgetful functor $U$ need not preserve limits, so that, in general, it need not be a right adjoint (see \cite[Section 3.6]{joint}). 
\par
Ben Wieland has found an example of a profinite group $G'$, with finite vcd and a closed normal subgroup $H'$, and 
a discrete $G'$-spectrum $Y$, such that $Y_{f,G'}$ is not fibrant as a discrete $H'$-spectrum. For, if 
$Y_{f,G'}$ is fibrant as a discrete $H'$-spectrum, then the map $\psi(Y)_{H'}^{G'}$ of (\ref{sigh}) is a weak 
equivalence. But $\psi(Y)_{H'}^{G'}$ is not a weak equivalence. The details of this example are 
in the Appendix. This example proves that, when $H$ is closed and not open in $K$, a fibrant 
object in $\mathrm{Spt}_ K$ is not necessarily fibrant in $\mathrm{Spt}_H$.
\par
The above considerations show that, in general, there is no reason to expect that $X_{f,K}$ 
is fibrant in $\mathrm{Spt}_H$, and, without this fact, the above proof 
of Theorem \ref{finite} fails to work.
\par
Now suppose that $K$ has finite vcd. By \cite[Theorem 4.2]{fibrantmodel}, there is a map 
\[X \rightarrow \colim_{N \vartriangleleft_o H} (\holim_\Delta \Gamma_K^\bullet (X_{f,K}))^N\] 
that is a weak equivalence in $\mathrm{Spt}_H$, such that the target is a fibrant discrete $H$-spectrum. 
Then, by 
\cite[Corollary 5.4]{fibrantmodel}, 
we can make the 
identification 
\[X^{hH} = (\colim_{N \vartriangleleft_o H} (\holim_\Delta \Gamma_K^\bullet (X_{f,K}))^N)^H \cong 
(\holim_\Delta \Gamma^\bullet_K(X_{f,K}))^H,\] so that $X^{hH}$ is 
the $K/H$-spectrum 
$(\holim_\Delta \Gamma^\bullet_K(X_{f,K}))^H$. 
\par
Notice that $K/H$ is an infinite profinite group: if not, then 
$H$ is a closed subgroup of $K$ with finite index, so $H$ must be open in $K$, contrary to our 
assumption. To form the $K/H$-homotopy fixed points of $X^{hH}$, we need to know that 
$X^{hH}$ can be regarded as either a discrete $K/H$-spectrum or as a homotopy limit of a 
tower of discrete $K/H$-spectra (or, more generally, as in \cite{joint}, the 
homotopy limit of a cofiltered diagram of discrete $K/H$-spectra), for these are the cases 
where the $K/H$-homotopy fixed points are defined. However, in general, it is not known 
how to do this. For example, the Appendix shows that $Y^{hH'}$ can not be regarded as a 
discrete $G'/H'$-spectrum (if it could be, then $\pi_\ast(Y^{hH'})$ is a discrete $G'/H'$-module, which 
is false), 
so that, in general, it can not be assumed that $X^{hH}$ is a discrete $K/H$-spectrum. Also, it 
is not known, in general, how to realize $X^{hH}$ as the homotopy limit of a tower or cofiltered 
diagram of discrete $K/H$-spectra. (For example, letting 
$\prod^\ast$ denote cosimplicial replacement,
\begin{align*}
X^{hH} & \cong \holim_\Delta (\Gamma^\bullet_K(X_{f,K}))^H = \mathrm{Tot}(\textstyle{\prod}^\ast 
\displaystyle(\Gamma_K^\bullet(X_{f,K}))^H) \\ 
& \cong \lim_n \mathrm{Tot}_n(\textstyle{\prod}^\ast \displaystyle(\Gamma_K^\bullet(X_{f,K}))^H)
\end{align*} presents $X^{hH}$ as the limit of a tower of $K/H$-spectra, with each 
$(\Gamma^l_K(X_{f,K}))^H$ a discrete $K/H$-spectrum (by the discussion just above Lemma 
\ref{fibrant}). However, 
$\mathrm{Tot}_n(\textstyle{\prod}^\ast (\Gamma_K^\bullet(X_{f,K}))^H)$ is not necessarily a 
discrete $K/H$-spectrum, since the infinite product $\prod^l(\Gamma_K^\bullet(X_{f,K}))^H$ need 
not be a discrete $K/H$-spectrum, so this attempt to present $X^{hH}$ as a continuous $K/H$-spectrum 
fails.)
\par
Let us consider in more detail whether or not $X^{hH}$ can be regarded as a discrete $K/H$-spectrum. 
Note that 
\[X^{hH} \cong \holim_\Delta (\Gamma^\bullet_K(X_{f,K}))^H \cong \holim_\Delta 
\colim_{N \vartriangleleft_o K} (\Gamma^\bullet_K(X_{f,K}))^{NH}.\] Also, there is a canonical 
$K/H$-equivariant map 
\[\Psi(X)_H^K \: \colim_{N \vartriangleleft_o K} X^{hNH} \rightarrow X^{hH}\] that is defined to be 
\[\colim_{N \vartriangleleft_o K} X^{hNH} \cong 
\colim_{N \vartriangleleft_o K} \holim_\Delta (\Gamma^\bullet_K(X_{f,K}))^{NH} 
\rightarrow \holim_\Delta \colim_{N \vartriangleleft_o K} (\Gamma^\bullet_K(X_{f,K}))^{NH}.\] 
Notice that the canonical map $K/H \rightarrow K/(NH)$ makes the source of 
$\Psi(X)_H^K$ a discrete 
$K/H$-spectrum, since $\holim_\Delta (\Gamma^\bullet_K(X_{f,K}))^{NH}$ is a 
$K/{(NH)}$-spectrum and $K/{(NH)}$ is a finite group. 
\par
If $\Psi(X)_H^K$ is a weak equivalence, then $X^{hH}$ can be identified with the discrete 
$K/H$-spectrum $\colim_{N \vartriangleleft_o K} \holim_\Delta (\Gamma^\bullet_K(X_{f,K}))^{NH}$, 
and, thus, in this case, $(X^{hH})^{hK/H}$ is defined (as 
$(\colim_{N \vartriangleleft_o K} \holim_\Delta (\Gamma^\bullet_K(X_{f,K}))^{NH})^{hK/H}$). 
However, once again, $\Psi(X)_H^K$ does not have to be a weak equivalence: if $\Psi(Y)_{H'}^{G'}$ 
were a weak equivalence, then $\pi_0(Y^{hH'})$ must be a discrete $G'/H'$-module, and 
this is false. 
\par
Another way to express the failure of $\Psi(X)_H^K$ to always be a weak equivalence is by 
saying that the above holim and filtered colimit do not have to commute. However, there are 
conditions that guarantee that the holim and colimit do commute (up to weak equivalence). Let 
$J$ be any closed subgroup of $K$. Then there is a homotopy spectral sequence $E_r^{\ast, \ast}(J)$ 
that has the form 
\[E_2^{s,t}(J) = H^s_c(J; \pi_t(X)) \Rightarrow \pi_{t-s}(\holim_\Delta (\Gamma^\bullet_K(X_{f,K}))^J)\] 
(see the proof of \cite[Lemma 7.12]{cts}). 
Notice that there is a map of spectral sequences 
\[\colim_{N \vartriangleleft_o K} E_r^{\ast, \ast}(NH) \rightarrow E_r^{\ast, \ast}(H),\] 
such that
\[\colim_{N \vartriangleleft_o K} E_2^{s,t}(NH) = \colim_{N \vartriangleleft_o K} H^s_c(NH; \pi_t(X)) 
\cong H^s_c(H; \pi_t(X)) = E_2^{s,t}(H).\] 
\par
By \cite[Proposition 3.3]{Mitchell}, if  
\begin{itemize}
\item[(i)]
there exists 
a fixed $p$ such that $H^s_c(NH; \pi_t(X)) = 0$ for all $s > p$, all $t \in \mathbb{Z}$, and all 
$N \negmedspace \vartriangleleft_o \negmedspace K$; or 
\item[(ii)]
there exists a fixed $q$ such that $H^s_c(NH; \pi_t(X)) = 0$ for 
all $t  > q$, all $s \geq 0$, and all $N$,
\end{itemize} then the colimit of spectral sequences has abutment equal to the colimit of the 
abutments, so that $\Psi(X)_H^K$ is a weak equivalence, and, hence, $X^{hH}$
can, as above, be regarded as a discrete $K/H$-spectrum.
\par
The two conditions above imply that $X^{hH}$ can be identified with the discrete $K/H$-spectrum 
$\colim_{N \vartriangleleft_o K} X^{hNH}$ if, for example, one of the following conditions is satisfied: 
\begin{itemize}
\item
$K$ has finite cohomological dimension (that is, there exists some $p$ such that 
$H^s_c(K; M) = 0$ for all $s >p$ and all discrete $K$-modules $M$), so that 
there is a uniform bound on 
the cohomological dimension of all the $NH$; or
\item
there exists some $q$ such that $\pi_t(X) = 0$ for all $t > q$.
\end{itemize}
We conclude that when $K$ has finite vcd, unless certain additional 
hypotheses are in place, it is not known 
how to show that $X^{hH}$ is a continuous $K/H$-spectrum, so that, 
in such situations, it is not known how to form the iterated homotopy fixed point 
spectrum $(X^{hH})^{hK/H}$, with respect to the natural $K/H$-action.
\par
The above discussion used the hypercohomology spectra $\holim_\Delta 
(\Gamma_K^\bullet (X_{f,K}))^J$, for various $J$ closed in $K$. But suppose 
we are in a situation where these are not available to us. (If $K$ has finite vcd, these are 
always available. If $K$ does not have finite vcd, there are other conditions that permit the use of the 
hypercohomology spectra: these conditions are similar to those listed above and 
can be deduced from considering the proofs of Theorem 7.4 and Lemma 7.12 in \cite{cts}.) 
Further, suppose that 
$X_{f,K}$ is not fibrant in $\mathrm{Spt}_H$ 
and that $(X_{f,K})^H \not\simeq (X_{f,H})^H$, so that the discrete $K/H$-spectrum 
$(X_{f,K})^H = \colim_{N \vartriangleleft_o K} (X_{f,K})^{NH} = 
\colim_{N \vartriangleleft_o K} X^{hNH}$ fails to be a model for $X^{hH}$. 
\par
Under these assumptions, the only model for $X^{hH}$ that is available to us 
is the definition $(X_{f,H})^H$, and, in this case, the abstract $X_{f,H}$ is not known 
to carry a $K$-action, so that $(X_{f,H})^H$ is not known to 
be a $K/H$-spectrum (with non-trivial $K/H$-action). In such a situation, apart from assigning 
$X^{hH}$ the trivial $K/H$-action (which typically is not a natural thing to do), 
it is not known how to consider its $K/H$-homotopy fixed points.
\section{Hyperfibrant and totally hyperfibrant discrete $G$-spectra}\label{hyperfibrant}
\par
In this section, 
we consider a situation in which the iterated homotopy fixed point 
spectrum is always defined. As before, 
throughout this section, we let $H$ and $K$ be closed subgroups 
of the profinite group $G$, with $H$ normal in $K$. 
\par
It is easy 
to see that the map $\psi \: \colim_{N} X^{hNK} \rightarrow 
X^{hK}$ of (\ref{sigh}) is a weak equivalence for all $K$ open in $G$, since 
$X_{f,G}$ is always fibrant in $\mathrm{Spt}_K$, for such $K$. If 
$X_{f,G}$ is fibrant in $\mathrm{Spt}_K$ for all closed $K$, then $\psi$ 
is a weak equivalence for all closed $K$, but, in the Appendix, Wieland and the 
author prove that this does not always 
happen. This motivates the following definition.
\begin{Def}\label{hyper}
Let $X$ be a discrete $G$-spectrum. If $\psi$ is a weak 
equivalence for all $K$ closed in $G$, then $X$ is a {\em hyperfibrant} 
discrete $G$-spectrum. 
\end{Def}
\par
As alluded to above, there do exist profinite groups $G$ for which there is a discrete 
$G$-spectrum that is not a hyperfibrant discrete $G$-spectrum: the map 
$\psi(\bigvee_{n \geq 0} \Sigma^n H(\gpring))_{\mathbb{Z}/p}^{\mathbb{Z}/p \times \mathbb{Z}_q}$, 
considered in the Appendix, is not a weak equivalence. 
\par
For the example below and the paragraph after it, we assume that $G$ has 
finite vcd.
\begin{Ex}\label{map}
Let $X$ be any spectrum (with no $G$-action), so that $\mathrm{Map}_c(G,X)$, 
as defined in Section \ref{summary}, is a discrete $G$-spectrum. The descent spectral sequence for 
$\mathrm{Map}_c(G,X)^{hK}$ has the form 
\[H^s_c(K;\mathrm{Map}_c(G,\pi_t(X))) \Rightarrow 
\pi_{t-s}(\mathrm{Map}_c(G,X)^{hK}).\] Since 
$H^s_c(K;\mathrm{Map}_c(G,\pi_t(X))) = 0$, for $s > 0$, 
\begin{align*}
\pi_\ast(\mathrm{Map}_c(G,X)^{hK}) & \cong 
H^0_c(K;\mathrm{Map}_c(G,\pi_\ast(X))) = \mathrm{Map}_c(G,\pi_\ast(X))^K \\ &
\cong \mathrm{Map}_c(G/K, \pi_\ast(X)) \cong 
\pi_\ast(\mathrm{Map}_c(G,X)^K),
\end{align*}
and thus, $\mathrm{Map}_c(G,X)^{hK} \simeq \mathrm{Map}_c(G,X)^K,$ for any 
$K$ closed in $G$. Therefore, we have:
\begin{align*}
\mathrm{Map}_c(G,X)^{hK} & \simeq \mathrm{Map}_c(G,X)^K \cong 
\colim_N \mathrm{Map}_c(G,X)^{NK} \\ & \simeq \colim_N 
\mathrm{Map}_c(G,X)^{hNK}.
\end{align*}
Thus, for any spectrum $X$, $\mathrm{Map}_c(G,X)$ is a hyperfibrant 
discrete $G$-spectrum.
\end{Ex}
\par
For any discrete $G$-spectrum $X$, 
\[X^{hK} \cong \holim_\Delta \colim_{N} 
(\Gamma^\bullet_G(X_{f,G}))^{NK},\] and similarly, 
\[\colim_{N} X^{hNK} = 
\colim_{N} \holim_\Delta 
(\Gamma^\bullet_G(X_{f,G}))^{NK}.\] Thus, if $\psi$ is a weak equivalence, 
for some $K$, then the canonical map (by a slight abuse of notation, we also call this map $\psi$)
\begin{equation}\label{anotherpsi}\zig
\psi \: \colim_{N} \holim_\Delta 
(\Gamma^\bullet_G(X_{f,G}))^{NK} \rightarrow \holim_\Delta 
\colim_{N} (\Gamma^\bullet_G(X_{f,G}))^{NK}
\end{equation}
is a weak equivalence. Hence, if $X$ is hyperfibrant, then, for all $K$, 
the colimit and the holim involved in defining $\psi$ commute.  
\begin{Lem}\label{colimit}
If $X$ is a discrete $G$-spectrum and $K$ is normal in $L$, where $L$ is a 
closed subgroup of $G$, then 
$\colim_{N} X^{hNK}$ is a discrete 
$L/K$-spectrum.
\end{Lem}
\begin{proof}
Note that the spectrum $X^{hNK} = (X_{f,G})^{NK}$ is an 
$NL/NK$-spectrum, where $NL/NK$ is a finite group. The continuous map 
$L/K \rightarrow NL/NK$, given by $lK \mapsto lNK$, 
makes $(X_{f,G})^{NK}$ a discrete $L/K$-spectrum. Since the forgetful 
functor $\mathrm{Spt}_{L/K} \rightarrow \mathrm{Spt}$ is a left adjoint, 
colimits in $\mathrm{Spt}_{L/K}$ are formed in the category of spectra, so that 
$\colim_{N} (X_{f,G})^{NK}$ is a discrete 
$L/K$-spectrum. 
\end{proof}
\par
Given a hyperfibrant discrete $G$-spectrum, there is a natural way to define the 
$K/H$-homotopy fixed points of $X^{hH}$.
\begin{Def}\label{construction}
Let $X$ be a hyperfibrant discrete $G$-spectrum, so that 
the map $\psi(X)_H^G \: \colim_N X^{hNH} \rightarrow X^{hH}$ is a weak equivalence. Notice 
that the source of this map, $\colim_N X^{hNH}$, is a discrete $K/H$-spectrum, and 
the closely related $K/H$-equivariant map 
$\psi$ of (\ref{anotherpsi}), also from $\colim_N X^{hNH}$ to $X^{hH}$, is a weak equivalence. 
Thus, we can identify $X^{hH}$ with the discrete $K/H$-spectrum 
$\colim_{N} X^{hNH},$ and, hence, the $K/H$-homotopy fixed points of $X^{hH}$ are given by 
\[(X^{hH})^{hK/H} = (\colim_{N} X^{hNH})^{hK/H}.\]
\end{Def}
\par
Let $X$ be any discrete $K$-spectrum. As in the proof of Lemma \ref{colimit}, $\colim_N X^{NH},$ where the colimit is over all open normal subgroups of $K$, is a discrete $K/H$-spectrum, so that the 
isomorphism $X^{H} \cong \colim_N X^{NH}$ implies that $X^H$ is also a 
discrete $K/H$-spectrum. Hence, there is a functor 
$(-)^H \: \mathrm{Spt}_{K} \rightarrow \mathrm{Spt}_{K/H}.$  
\par
The author thanks Mark Behrens for help with the next lemma, which is basically a version of \cite[Lemma 6.35]{Jardine}, which is for 
simplicial presheaves.
\begin{Lem}\label{fibrant}
The functor $(-)^H \: \mathrm{Spt}_{K} \rightarrow \mathrm{Spt}_{K/H}$ 
preserves fibrant objects.
\end{Lem}
\begin{proof}
It is easy to see that 
the functor $(-)^H$ has a left adjoint $t \: \mathrm{Spt}_{K/H} 
\rightarrow \mathrm{Spt}_K$ that sends a discrete $K/H$-spectrum $X$ to 
$X$, where now $X$ is regarded as a discrete $K$-spectrum via the canonical map 
$K \rightarrow K/H$. It suffices to show that $t$ preserves all weak 
equivalences and cofibrations, since this implies that $(-)^H$ preserves 
all fibrations. 
\par
If $f$ is a weak equivalence (cofibration) in $\mathrm{Spt}_{K/H}$, then 
$f$ is a $K/H$-equivariant map that is a weak equivalence (cofibration) of 
spectra. Since, as a map of spectra, $t(f) = f$, $t(f)$ is a weak equivalence (cofibration) of spectra. Because $t(f)$ is also $K$-equivariant, it is a weak equivalence 
(cofibration) in $\mathrm{Spt}_K$.
\end{proof}
\par
This lemma implies the following result, giving another useful property of 
hyperfibrant discrete $G$-spectra.
\begin{Lem}\label{easycase}
If $X$ is a hyperfibrant discrete $G$-spectrum and $K$ is any closed normal subgroup of 
$G$, then $(X^{hK})^{hG/K} \simeq X^{hG}$.
\end{Lem}
\begin{proof}
This result follows immediately from Lemma \ref{iterated} below, by letting 
$K = G$ and $H = K$.
\end{proof}
\par
Now let $X$ be a hyperfibrant discrete $G$-spectrum and suppose that 
the closed subgroup $K$ is a proper subgroup of $G$. Let us consider if there is 
an equivalence $(X^{hH})^{hK/H} \simeq X^{hK}$. Notice that
\[(X^{hH})^{hK/H} = (\colim_N (X_{f,G})^{NH})^{hK/H} = 
(((X_{f,G})^H)_{f,K/H})^{K/H}\] and 
\[X^{hK} \overset{\simeq}{\longleftarrow} \colim_N (X_{f,G})^{NK} = (X_{f,G})^K = ((X_{f,G})^H)^{K/H}.\]
Thus, $(X^{hH})^{hK/H} \simeq X^{hK}$ if the canonical 
map 
\begin{equation}\label{induced}\zig
((X_{f,G})^H)^{K/H} \rightarrow (((X_{f,G})^H)_{f,K/H})^{K/H}
\end{equation} is 
a weak equivalence. 
\par
If $(X_{f,G})^H$ were fibrant in $\mathrm{Spt}_{K/H}$, then, in $\mathrm{Spt}_{K/H}$, 
the map 
\[(X_{f,G})^H \rightarrow ((X_{f,G})^H)_{f,K/H}\] is a weak equivalence between fibrant objects, so that the map in (\ref{induced}) is a weak equivalence (by \cite[Corollary 3.9]{cts}). 
Thus, if $X_{f,G}$ is fibrant in $\mathrm{Spt}_K$, then the preceding observation and Lemma 
\ref{fibrant} imply that $(X^{hH})^{hK/H} \simeq X^{hK}$. But, as we saw with 
$Y_{f,G'}$ not being fibrant in $\mathrm{Spt}_{H'}$, the forgetful 
functor $\mathrm{Spt}_G \rightarrow \mathrm{Spt}_K$ does not necessarily preserve fibrant objects. 
Also, we do not know of any general argument, for $K$ not open in $G$, that would show that 
$(X_{f,G})^H$ is fibrant in $\mathrm{Spt}_{K/H}$. In conclusion, though we do not have an 
example of $(X^{hH})^{hK/H} \simeq X^{hK}$ failing to hold when $X$ is a hyperfibrant discrete 
$G$-spectrum, we also do not know of any general argument that shows that, when 
$X$ is hyperfibrant, this equivalence must always hold. This leads us to the 
following considerations.
\begin{Lem}\label{iterated}
Let $X$ be a discrete $G$-spectrum and let $K$ be a closed subgroup of $G$. If $X$ is a hyperfibrant discrete $K$-spectrum, then there is a weak 
equivalence $(X^{hH})^{hK/H} \simeq X^{hK}$, where $H$ is any closed subgroup of 
$G$ that is also normal in $K$.
\end{Lem}
\begin{proof}
Since $X$ is hyperfibrant in $\mathrm{Spt}_K$, $\psi \: \colim_N X^{hNH} \rightarrow 
X^{hH}$ is a weak equivalence, where the colimit is over all open normal 
subgroups of $K$. Thus, by Definition \ref{construction},
\[(X^{hH})^{hK/H} = (\colim_N X^{hNH})^{hK/H} = 
(\colim_N (X_{f,K})^{NH})^{hK/H} \simeq ((X_{f,K})^H)^{hK/H}.\]
Since $X_{f,K}$ is fibrant in $\mathrm{Spt}_K$, Lemma \ref{fibrant} 
implies that $(X_{f,K})^H$ is fibrant in $\mathrm{Spt}_{K/H}$. Also, the identity map 
$(X_{f,K})^H \rightarrow (X_{f,K})^H$ is a trivial cofibration in $\mathrm{Spt}_{K/H}$. Thus, we can 
let $(X_{f,K})^H = ((X_{f,K})^H)_{f,K/H},$ and hence, 
\[((X_{f,K})^H)^{hK/H} = (((X_{f,K})^H)_{f,K/H})^{K/H} = 
((X_{f,K})^H)^{K/H} = (X_{f,K})^K = X^{hK}.\]
\end{proof}
\par
This lemma implies that if the discrete $G$-spectrum $X$ 
is hyperfibrant as a discrete $K$-spectrum, for all $K$ closed in $G$, then 
not only is the iterated homotopy fixed point spectrum always defined, 
but (\ref{question}) is always true. Thus, we make the following definition. 
\begin{Def}\label{totallyhyper}
Let $X$ be a discrete $G$-spectrum. If $X$ is a hyperfibrant discrete 
$K$-spectrum, for all $K$ closed in $G$, then $X$ is a 
{\em totally hyperfibrant} discrete $G$-spectrum. Therefore, if $X$ is a 
totally hyperfibrant discrete $G$-spectrum, then $(X^{hH})^{hK/H} \simeq 
X^{hK}$, for any $H$ and $K$ closed in $G$, with $H$ normal in $K$.
\end{Def}
\begin{Ex}\label{totally}
For any spectrum $X$, $\mathrm{Map}_c(G,X)$ is a totally hyperfibrant 
discrete $G$-spectrum. To see this, let 
$K$ be any closed subgroup of $G$, and let $L$ be 
any closed subgroup of $K$. 
Since $L$ is closed in $G$, $\mathrm{Map}_c(G,X)^{hL}$ is defined. Also, if $N$ is an open normal subgroup of $K$, then $NL$ is closed in $K$, 
since $K$ is a profinite group and $NL$ is an open subgroup of $K$. Thus, 
$NL$ is a closed subgroup of $G$, so that $\mathrm{Map}_c(G,X)^{hNL}$ is defined. Then, by Example \ref{map} and by taking the colimit over all open normal 
subgroups $N$ of $K$, we have: 
\[ \colim_N \mathrm{Map}_c(G,X)^{hNL} \simeq \colim_N \mathrm{Map}_c(G,X)^{NL} 
\cong \mathrm{Map}_c(G,X)^L \simeq \mathrm{Map}_c(G,X)^{hL}.\]
\end{Ex}
\section{Iterated homotopy fixed points for $E_n$}
In this section, we show that it is always possible to construct the 
spectrum $(E_n^{dhH})^{hK/H}$. As mentioned in the Introduction, applying 
(\ref{thesis}) allows us to identify $E_n^{dhH}$ with 
$E_n^{hH}$, yielding the iterated homotopy fixed point spectrum 
$(E_n^{hH})^{hK/H}$. We begin with some notation. 
\par
We use $(-)_\mathtt{f}$ to 
denote functorial fibrant replacement in the category of spectra. Let $E(n)$ be the 
Johnson-Wilson spectrum, with $E(n)_\ast = 
\mathbb{Z}_{(p)}[v_1, ..., v_{n-1}][v_n^{\pm 1}],$ where the degree of each $v_i$ is $2(p^i-1)$.   Let 
\[M_{0} \leftarrow M_{1} \leftarrow M_{2} \leftarrow \cdots\]
be a tower of generalized Moore spectra such that $\lhat(S^0) 
\simeq \holim_i (L_n M_{i})_\mathtt{f}$, where $L_n$ is Bousfield localization with respect to $E(n)$. This tower of finite type $n$ spectra exists by \cite[Proposition 
4.22]{HS}. Let $\{U_j\}_{j \geq 0}$ be a descending chain of open normal 
subgroups of $G_n,$ with $\bigcap_j U_j = \{e\}$ (see \cite[(1.4)]{DH}), so 
that $G_n \cong \lim_j G_n/U_j.$ Also, 
let $H$ and $K$ be closed subgroups of $G_n$, with $H$ normal in $K$. 
\par
We recall the following useful 
result from {\cite[Section 2]{HoveyCech}}.
\begin{Lem}\label{complete}
If $Z$ is an $E(n)$-local spectrum, then there is an equivalence
\[\lhat Z \simeq \holim_i(Z \wedge M_{i})_\mathtt{f}.\]
\end{Lem}
\par
By \cite[Definition 1.5]{DH}, $E_n^{dhH} \simeq \lhat(\colim_j E_n^{dhU_jH})$, and, by 
\cite[Lemma 7.2]{HS}, there is an equivalence 
$\lhat(\colim_j E_n^{dhU_jH}) \wedge 
M_{i} \simeq L_n(\colim_j E_n^{dhU_jH}) \wedge M_{i}.$ Thus, since 
each $E_n^{dhU_jH}$ is $E(n)$-local, we immediately obtain 
the following result. 
\begin{Lem}\label{hyper2}
For each $i$, $E_n^{dhH} \wedge M_{i} \simeq 
\colim_j (E_n^{dhU_jH} \wedge M_{i}).$
\end{Lem}
\par
The following definition, which we recall from \cite{cts}, defines a useful spectrum.
\begin{Def}
Let $F_n = \colim_j E_n^{dhU_j}$, a discrete $G_n$-spectrum. 
\end{Def}
\begin{Rk}
Using \cite[Corollary 9.8]{cts} and (\ref{thesis}), Lemma \ref{hyper2} implies that
\begin{align*}
(F_n \wedge M_{i})^{hH} & \simeq E_n^{dhH} \wedge M_{i} \simeq 
\colim_j (E_n^{dhU_jH} \wedge M_{i}) \\ & \simeq 
\colim_j (F_n \wedge M_{i})^{hU_jH} \simeq 
\colim_{N} (F_n \wedge M_{i})^{hNH},
\end{align*}
where the last colimit is over all open normal subgroups of $G_n$. 
Therefore, since Lemma \ref{hyper2} is valid if $H$ is 
replaced with any closed subgroup of $G_n$, $F_n \wedge M_{i}$ 
is a hyperfibrant discrete $G_n$-spectrum.
\end{Rk}
\par
For the rest of this section, we let $X$ be any spectrum, with no $K/H$-action.
\begin{Lem}
The spectrum $\colim_j (E_n^{dhU_jH} \wedge M_{i} \wedge X)_\mathtt{f}$ 
is a discrete $K/H$-spectrum. 
\end{Lem}
\begin{proof}
Note that $E_n^{dhU_jH}$ is a $U_jK/U_jH$-spectrum. Thus, 
$E_n^{dhU_jH} \wedge M_{i} \wedge X$ is a $U_jK/U_jH$-spectrum. By 
functoriality, $(E_n^{dhU_jH} \wedge M_{i} \wedge X)_\mathtt{f}$ is also a 
$U_jK/U_jH$-spectrum. The argument is completed as in the proof of Lemma 
\ref{colimit}.
\end{proof} 
\begin{Cor}\label{cor}
The spectrum
\[\lhat(E_n^{dhH} \wedge X) \simeq 
\holim_i  (\colim_j (E_n^{dhU_jH} \wedge M_{i} \wedge X)_\mathtt{f})\] is a 
continuous $K/H$-spectrum. In particular, \[E_n^{dhH} \simeq 
\holim_i(\colim_j(E_n^{dhU_jH} \wedge M_{i})_\mathtt{f})\] is a continuous 
$K/H$-spectrum. 
\end{Cor}
\begin{proof}
By Lemma \ref{complete},
\begin{align*}
\lhat(E_n^{dhH} \wedge X) & \simeq \holim_i 
(E_n^{dhH} \wedge M_{i} \wedge X)_\mathtt{f} \\ 
& \simeq \holim_i (\colim_j (E_n^{dhU_jH} \wedge M_{i} \wedge X)_\mathtt{f}).
\end{align*} 
Since $(E_n^{dhU_jH} \wedge M_{i} \wedge X)_\mathtt{f}$ is fibrant, so is 
$\colim_j (E_n^{dhU_jH} \wedge M_{i} \wedge X)_\mathtt{f}$, so that the homotopy limit is of a diagram of fibrant spectra, as required.
\end{proof}
\begin{Rk}
We can just as well phrase Corollary \ref{cor} as saying that 
\[\lhat(E_n^{dhH} \wedge X) \simeq 
\holim_i \, 
((\colim_j E_n^{dhU_jH}) \wedge M_{i} \wedge X)_{f,K/H}.\] In 
particular, note that $E_n^{dhH} \simeq 
\holim_i ((\colim_j E_n^{dhU_jH}) \wedge M_{i})_{f,K/H}.$ Now recall 
that $E_n \simeq \holim_i (F_n \wedge M_{i})_{f,G_n}.$ 
Thus, 
$\colim_j E_n^{dhU_jH}$ is playing the same role as $F_n$ and is the 
analogue of $F_n.$
\end{Rk}
\par
Note that, since $\lhat(E_n^{dhH} \wedge X)$ is a continuous $K/H$-spectrum (that is, 
we identify this spectrum with 
$\holim_i  (\colim_j (E_n^{dhU_jH} \wedge M_{i} \wedge X)_\mathtt{f})$), 
Definition \ref{hg} gives
\[(\lhat(E_n^{dhH} \wedge X))^{hK/H} = \holim_i((\colim_j (E_n^{dhU_jH} \wedge M_{i} \wedge X)_\mathtt{f})^{hK/H}).\]
\section{The descent spectral sequence for $\lhat(E_n^{dhH} \wedge X)$}\label{dsssection}
Let $H$ and $K$ be closed subgroups of $G_n$, with $H$ normal in $K$. 
By \cite[Theorem 9.6]{Dixon}, if $G$ is a $p$-adic analytic group, 
$U$ a closed subgroup of $G$, and $N$ a closed normal subgroup of $G$, 
then $U$ and $G/N$, with the quotient topology, are $p$-adic 
analytic groups. Thus, we obtain the following useful result.
\begin{Lem}\label{vcd}
The profinite group $K/H$ 
is a compact $p$-adic analytic group and $\mathrm{vcd}(K/H) < \infty.$ 
\end{Lem}
\par
Since, by Corollary \ref{cor}, $\lhat(E_n^{dhH} \wedge X)$ is a continuous 
$K/H$-spectrum, for any spectrum $X$ with trivial $K/H$-action, 
Theorem \ref{dss} gives a descent spectral sequence
\begin{equation}\zig\label{dss2}
E_2^{s,t} \Rightarrow \pi_{t-s}((\lhat(E_n^{dhH} \wedge X))^{hK/H}),
\end{equation}
where \[E_2^{s,t}=\pi^s\pi_t(\holim_i 
(\Gamma^\bullet_{K/H}(\colim_j (E_n^{dhU_jH} 
\wedge M_{i} \wedge X)_\mathtt{f}))^{K/H}).\]
\par
Let $H^s_\mathrm{cts}(G;M)$ denote the continuous cohomology of continuous cochains of $G$, with coefficients in the topological $G$-module $M$ 
(see \cite[Section 2]{Tate}). Given a tower $\{M_i\}$ of discrete $G$-modules, 
\cite[Theorem 2.2]{Jannsen} implies that, if the tower of abelian groups $\{M_i\}$ satisfies the 
Mittag-Leffler condition, then 
\[H^s_\mathrm{cont}(G ; \{M_i\}) \cong H^s_\mathrm{cts}(G; \lim_i M_i),\] for all $s \geq 0.$ Therefore, 
if the tower of abelian groups $\{\pi_t(E_n^{dhH} \wedge M_{i} \wedge X)\}_i$ satisfies the Mittag-Leffler condition for each $t \in \mathbb{Z}$, then 
\[E_2^{s,t} \cong 
H^s_\mathrm{cont}(K/H;\{\pi_t(\lhat(E_n^{dhH} \wedge M_{i} \wedge X))\}) \cong 
H^s_\mathrm{cts}(K/H;\pi_t(\lhat(E_n^{dhH} \wedge X))).\] 
\par
Now let $X$ be a finite spectrum. Then 
$\lhat(E_n^{dhH} \wedge X) \simeq E_n^{dhH} \wedge X$ and the hypotheses of the preceding 
sentence are satisfied, so that 
\[E_2^{s,t} \cong H^s_\mathrm{cts}(K/H;\pi_t(E_n^{dhH} \wedge X)).\] 
Additionally, $\pi_t(E_n^{dhH} \wedge X) \cong \lim_i \pi_t(E_n^{dhH} \wedge M_i \wedge X)$ is a 
profinite $\mathbb{Z}_p[[K/H]]$-module, and, in this case, the continuous cohomology has 
more structure. In general, if $M$ is a profinite $\mathbb{Z}_p[[K/H]]$-module with 
$M = \lim_\alpha M_\alpha$, where each $M_\alpha$ is a finite discrete 
$\mathbb{Z}_p[[K/H]]$-module, then 
\[H^s_\mathrm{cts}(K/H; M) \cong \lim_\alpha H^s_c(K/H; M_\alpha) =: H^s_c(K/H; M),\] 
for all $s \geq 0$, 
where the rightmost term can also be defined as the $\mathrm{Ext}$ group 
$\mathbf{Ext}_{\mathbb{Z}_p[[K/H]]}^s(\mathbb{Z}_p, M)$ 
(see \cite[pg. 137]{LHS} and \cite[(3.7.10)]{Symonds} for more detail). Putting the above facts together, 
we obtain that spectral sequence (\ref{dss2}) has the form 
\begin{equation}\zig\label{dss3}
H^s_c(K/H; \pi_t(E_n^{dhH} \wedge X)) \Rightarrow 
\pi_{t-s}((E_n^{dhH} \wedge X)^{hK/H}).
\end{equation}
\par
Recall from \cite{LHS} that there is a Lyndon-Hochschild-Serre 
spectral sequence
\begin{equation}\zig\label{lhs} 
H^s_c(K/H; \pi_t(E_n^{dhH} \wedge X)) \Rightarrow 
\pi_{t-s}(E_n^{dhK} \wedge X),
\end{equation}
where $X$ is a finite spectrum. 
The fact that spectral sequences (\ref{dss3}) and (\ref{lhs}) have identical 
$E_2$-terms suggests that there is an equivalence 
$(E_n^{dhH})^{hK/H} \simeq E_n^{dhK}.$ In the next result, we use 
descent spectral sequence (\ref{dss2}), when $X=E_n$, to show that 
this equivalence holds after 
smashing with $E_n$, before taking $K/H$-homotopy fixed 
points. First of all, we define the relevant map.
\par
Let \[C_i =\colim_j 
(E_n^{dhU_jH} \wedge M_{i} \wedge E_n)_\mathtt{f}.\]
Then there is a map 
\[\theta \: \lhat(E_n^{dhK} \wedge E_n) \rightarrow 
(\lhat(E_n^{dhH} \wedge E_n))^{hK/H},\]
which is defined by making the identifications
\[\lhat(E_n^{dhK} \wedge E_n) = \holim_i (\colim_j 
(E_n^{dhU_jK} \wedge M_{i} \wedge E_n)_\mathtt{f})\] and 
\[(\lhat(E_n^{dhH} \wedge E_n))^{hK/H} = \holim_i 
(((C_i)_{f,K/H})^{K/H}),\] and 
by taking the composite of the canonical map 
\[\holim_i (\colim_j 
(E_n^{dhU_jK} \wedge M_{i} \wedge E_n)_\mathtt{f}) \rightarrow 
\lim_{K/H} (\holim_i C_i),\] the map \[\lim_{K/H} (\holim_i C_i) \overset{\cong}{\longrightarrow} 
\holim_i \lim_{K/H} C_i \overset{\cong}{\longrightarrow} \holim_i ((C_i)^{K/H}),\] and the map 
\[\holim_i ((C_i)^{K/H}) \rightarrow 
\holim_i (((C_i)_{f,K/H})^{K/H}).\]
\par
To avoid confusion, we point out that, in the expression $(\lhat(E_n^{dhH} \wedge E_n))^{hK/H},$ the lone $E_n$ has trivial $K/H$-action, as stated just after Lemma 
\ref{vcd}.
\begin{Thm}
The map 
\begin{equation}\zig\label{iso}
\theta \: \lhat(E_n^{dhK} \wedge E_n) \overset{\simeq}{\longrightarrow} 
(\lhat(E_n^{dhH} \wedge E_n))^{hK/H}\end{equation} is a weak equivalence. 
\end{Thm}
\begin{proof}
By \cite[Proposition 6.3]{DH}, 
\[\pi_t(E_n^{dhH} \wedge E_n \wedge M_{i}) \cong \mathrm{Map}_c(G_n/H, 
\pi_t(E_n \wedge M_{i})).\] Thus, 
since $\{\pi_t(E_n \wedge M_{i})\}$ satisfies 
the Mittag-Leffler condition for every integer $t$, 
the tower $\{\pi_t(E_n^{dhH} \wedge E_n \wedge M_{i})\}$ 
does too, because $\mathrm{Map}_c(G_n/H,-) \: \mathbf{Ab} \rightarrow \mathbf{Ab}$, 
where $\mathbf{Ab}$ is the category of abelian groups, is an exact and additive functor. 
Therefore, descent spectral sequence (\ref{dss2}) takes the form 
\[E_2^{s,t} = 
H^s_\mathrm{cts}(K/H;\mathrm{Map}_c(G_n/H,\pi_t(E_n))) \Rightarrow 
\pi_{t-s}((\lhat(E_n^{dhH} \wedge E_n))^{hK/H}).\]
\par
Since the continuous epimorphism $\pi \: G_n/H \rightarrow G_n/K$ 
has a continuous section $\sigma \: G_n/K \rightarrow G_n/H$, such that $\sigma(eK) =eH$, 
there is a homeomorphism \[h \: K/H \times G_n/K \rightarrow 
G_n/H, \ \ \ (kH, gK) \mapsto kH \cdot \sigma(gK)\] (see \cite[proof of Lemma 3.15]{LHS}). 
Notice that $h$ is 
$K/H$-equivariant, with $K/H$ acting on the source by acting only on $K/H$, since 
\[k'H \cdot (kH,gK) = (k'kH,gK) \mapsto k'kH \cdot \sigma(gK) = k'H \cdot (kH \cdot \sigma(gK)).\] Thus, 
\[\mathrm{Map}_c(G_n/H,\pi_t(E_n)) \cong 
\mathrm{Map}_c(K/H,\mathrm{Map}_c(G_n/K, \pi_t(E_n)))\] is an 
isomorphism of topological $K/H$-modules. Therefore, 
\[E_2^{s,t} \cong \lim_i \, 
H^s_c(K/H; \mathrm{Map}_c(K/H,\mathrm{Map}_c(G_n/K, 
\pi_t(E_n \wedge M_{i})))),\] which is $0$, for $s>0$, and 
equals $\mathrm{Map}_c(G_n/K, \pi_t(E_n)),$ when $s=0$. 
Thus, \[\pi_\ast((\lhat(E_n^{dhH} \wedge E_n))^{hK/H}) 
\cong \mathrm{Map}_c(G_n/K,\pi_\ast(E_n)) \cong 
\pi_\ast(\lhat(E_n^{dhK} \wedge E_n)).\]
\end{proof}
\begin{Rk}
In Theorem \ref{main2}, by using (\ref{thesis}), 
we will show that 
$(E_n^{hH})^{hK/H} \simeq E_n^{hK}.$ Therefore, using (\ref{thesis}) again, 
(\ref{iso}) implies 
that taking $K/H$-homotopy fixed points commutes with smashing with 
$E_n$: 
\[(\lhat(E_n^{hH} \wedge E_n))^{hK/H} \simeq \lhat(E_n^{hK} \wedge E_n).\] 
This is 
interesting because such a commutation need not hold in general, unless one 
is smashing with a finite spectrum. However, (\ref{iso}) is not 
surprising, because it is known that, for all $H$,
\[\lhat(E_n^{dhH} \wedge E_n) \simeq (\lhat(E_n \wedge E_n))^{hH},\] 
where on the right-hand side, the second $E_n$ has the trivial $H$-action. This last 
equivalence follows from the fact that 
\[\pi_\ast(\lhat(E_n^{dhH} \wedge E_n)) \cong \mathrm{Map}_c(G_n, 
E_{n \ast})^H \cong \pi_\ast((\lhat(E_n \wedge E_n))^{hH}),\] where 
the second 
isomorphism is obtained from the descent spectral sequence for 
$(\lhat(E_n \wedge E_n))^{hH}$ and the fact that 
$\pi_\ast(\lhat(E_n \wedge E_n)) \cong \mathrm{Map}_c(G_n, E_{n \ast})$ (see 
the last paragraph of \cite{cts}).
\end{Rk}
\section{A proof of (\ref{question}), in the case of $E_n$, and some consequences}\label{morava}
\par
In \cite{thesis} (see \cite[Theorem 8.2.1]{joint} for a more efficient proof), the author showed that, as stated in (\ref{thesis}), $E_n^{hK}$ can be 
identified with $E_n^{dhK},$ for all closed subgroups $K$ of $G_n$. In this section, we use this result to show that (\ref{question}) holds when $X$ equals $E_n$.
\begin{Lem}\label{DH}
Let $K$ be any closed subgroup of $G_n$ and let $L$ be any closed subgroup of 
$K$. Then the canonical map $\colim_N (E_n^{dhNL} \wedge M_{i}) 
\rightarrow E_n^{dhL} \wedge M_{i}$, where the colimit is over all open 
normal subgroups $N$ of $K$, is a weak equivalence.
\end{Lem}
\begin{proof}
Both $L$ and $NL$ are closed subgroups of $G_n$, 
so that $E_n^{dhL}$ and $E_n^{dhNL}$ are defined. 
Since both the source and target of the map are $E_n$-local, it suffices to 
show that the associated map 
\[\colim_N (E_n^{dhNL} \wedge E_n \wedge M_{i}) 
\rightarrow E_n^{dhL} \wedge E_n \wedge M_{i}\] is a weak equivalence. For any integer $t$, by using \cite[Proposition 6.3]{DH} and the fact that 
$\mathrm{Map}_c(G_n,\pi_t(E_n \wedge M_{i}))$ is a discrete 
$K$-module (since it is a discrete $G_n$-module), we have: 
\begin{align*}
\pi_t(\colim_N (E_n^{dhNL} \wedge E_n \wedge M_{i})) & 
\cong \colim_N \mathrm{Map}_c(G_n, \pi_t(E_n \wedge 
M_{i}))^{NL} \\ & \cong \mathrm{Map}_c(G_n,  
\pi_t(E_n \wedge M_{i}))^L \\  
& \cong \pi_t(E_n^{dhL} \wedge E_n \wedge M_{i}),
\end{align*} completing the proof.
\end{proof} 
\begin{Cor}\label{th}
For each $i$, the spectrum $F_n \wedge M_{i}$ is a totally hyperfibrant discrete $G_n$-spectrum.
\end{Cor}
\begin{proof}
Let $K$ be any closed subgroup of $G_n$ and let $L$ be any closed subgroup 
of $K$. Then, by \cite[Corollary 9.8]{cts} and Lemma \ref{DH},
\begin{align*}
(F_n \wedge M_{i})^{hL} & \simeq E_n^{hL} \wedge M_{i} 
\simeq E_n^{dhL} \wedge M_{i} \simeq
\colim_N (E_n^{dhNL} \wedge M_{i}) 
\\ & \simeq \colim_N (E_n^{hNL} \wedge M_{i}) 
\simeq \colim_N (F_n \wedge M_{i})^{hNL},
\end{align*} 
where the colimit is over all 
open normal subgroups $N$ of $K$. 
\end{proof}
\par
Now let $H$ be any closed subgroup of $G_n$ that is also normal in $K$. 
Since $E_n^{hH} \simeq E_n^{dhH}$, we identify these two spectra, so 
that the construction of $(E_n^{dhH})^{hK/H}$ automatically yields the 
iterated homotopy fixed point spectrum $(E_n^{hH})^{hK/H}.$
\begin{Thm}\label{main2}
For any $H$ and $K$, there is an equivalence \[(E_n^{hH})^{hK/H} \simeq 
E_n^{hK}.\]
\end{Thm}
\begin{proof}
By applying Lemma \ref{iterated} to Corollary \ref{th}, we have that 
\[((F_n \wedge M_{i})^{hH})^{hK/H} \simeq (F_n \wedge M_{i})^{hK}.\]
Thus, we obtain: 
\begin{align*}
(E_n^{hH})^{hK/H} & = (\holim_i (F_n \wedge M_{i})^{hH})^{hK/H} 
= \holim_i ((F_n \wedge M_{i})^{hH})^{hK/H} \\ 
& \simeq \holim_i (F_n \wedge M_{i})^{hK} = E_n^{hK}.
\end{align*}
\end{proof}
\par
Let $X$ be any finite spectrum and continue to let $H$ and $K$ be closed subgroups of $G_n$, 
with $H \vartriangleleft K$. Also, recall from Lemma \ref{vcd} that, 
for any closed subgroups $J$ and $L$ in $G_n$, with $J \vartriangleleft L$, 
the profinite group $L/J$ has finite vcd. Then, by \cite[Remark 7.16]{cts}, $L/J$-homotopy 
fixed points commute with smashing with $X$. We use the equivalence $E_n^{hJ} \simeq 
E_n^{dhJ}$ as needed. Thus, (\ref{dss3}) gives a descent spectral 
sequence that has the form 
\begin{equation}\label{comp.one}\zig
H^s_c(K/H; \pi_t(E_n^{hH} \wedge X)) \Rightarrow \pi_{t-s}((E_n^{hH})^{hK/H} \wedge X).
\end{equation}
By \cite[(0.1)]{LHS} (see (\ref{lhs})), 
there is a strongly convergent Adams-type spectral sequence that has 
the form 
\begin{equation}\label{comp.two}\zig
H^s_c(K/H; \pi_t(E_n^{hH} \wedge X)) \Rightarrow \pi_{t-s}(E_n^{hK} \wedge X).
\end{equation} The following result shows that these two spectral sequences are isomorphic to 
each other.
\begin{Thm}\label{sscomparison}
The descent spectral sequence of $\mathrm{(}\negthinspace$\ref{comp.one}$\mathrm{)}$ is isomorphic to the strongly convergent 
spectral sequence of $\mathrm{(}$\negthinspace\ref{comp.two}$\mathrm{)}$ from the $E_2$-terms 
onward.
\end{Thm}
\begin{proof}
By Theorem \ref{main2}, the abutment of (\ref{comp.one}) can be written as 
$\pi_{t-s}(E_n^{hK} \wedge X)$. Then, it is not hard to see that spectral sequence 
(\ref{comp.one}) is an inverse limit over $i \geq 0$ of conditionally convergent 
descent spectral sequences 
$E_r^{\ast, \ast}(H, K, i)$ that 
have the form
\[E_2^{s,t}(H,K,i) = H^s_c(K/H; \pi_t(E_n^{hH} \wedge M_i \wedge X)) \Rightarrow 
\pi_{t-s}(E_n^{hK} \wedge M_i \wedge X)\] (e.g., see the proof of 
\cite[Proposition 7.4]{HMS}). Similarly, spectral sequence 
(\ref{comp.two}) is an inverse limit over $i \geq 0$ of strongly convergent 
Adams-type spectral sequences 
$\mathbb{E}_r^{\ast, \ast}(H,K,i)$ that have the same form as $E_r^{\ast, \ast}(H, K, i)$: 
\[\mathbb{E}_2^{s,t}(H,K,i) = H^s_c(K/H; \pi_t(E_n^{hH} \wedge M_i \wedge X)) \Rightarrow 
\pi_{t-s}(E_n^{hK} \wedge M_i \wedge X).\] Hence, to prove the theorem, it suffices 
to show that there is an isomorphism $E_r^{\ast, \ast}(H, K, i) \cong \mathbb{E}_r^{\ast, \ast}(H, K, i)$ 
of spectral sequences, from the $E_2$-terms onward, for each $i$.
\par
Notice that the limit $\lim_j U_jK/U_jH$ presents $K/H$ as a profinite group and there is an 
isomorphism 
\begin{align*} 
\colim_j \mathbb{E}_2^{s,t}(U_jH, U_jK, i) & = 
\colim_j H^s(U_jK/U_jH; \pi_t(E_n^{hU_jH} \wedge M_i \wedge X)) \\
& \cong H^s_c(K/H; \pi_t(E_n^{hH} \wedge M_i \wedge X)),
\end{align*} where the isomorphism applies Lemma \ref{hyper2}. Thus, noting that 
each spectral sequence $\mathbb{E}_r^{\ast,\ast}(U_jH, U_jK, i)$ has abutment 
$\pi_\ast(E_n^{hU_jK} \wedge M_i \wedge X)$ and using Lemma \ref{hyper2} again, it is easy to 
see that there is an isomorphism 
\begin{equation}\label{finaliso}\zig
\mathbb{E}_r^{\ast,\ast}(H, K, i) \cong \colim_j \mathbb{E}_r^{\ast,\ast}(U_jH, U_jK, i)
\end{equation} of 
spectral sequences. Now, since each $U_jK/U_jH$ is a finite group, \cite[Theorem A.1]{LHS} 
shows that there is an isomorphism 
\[\mathbb{E}_r^{\ast,\ast}(U_jH, U_jK, i) \cong E_r^{\ast,\ast}(U_jH, U_jK, i)\] of spectral sequences, 
for all $j$. Thus, 
\[\mathbb{E}_r^{\ast,\ast}(H, K, i) \cong \colim_j E_r^{\ast,\ast}(U_jH, U_jK, i) \cong 
E_r^{\ast, \ast}(H,K,i),\] where the last isomorphism is verified in the same way that 
(\ref{finaliso}) is.
\end{proof}
\par
The results in this paper apply to other 
spectra that are closely related to $E_n$. Let $k$ be any finite field containing 
$\mathbb{F}_{p^n}$. By \cite[Section 7]{Pgg/Hop0}, for any height $n$ formal 
group law $\Gamma$ over $k$, there is a commutative $S$-algebra $E(k, \Gamma)$ with 
\[E(k, \Gamma)_\ast = \mathbb{W}(k)[[u_1, ..., u_{n-1}]][u^{\pm 1}],\] 
graded in the same way that $E_{n\ast}$ is. 
Ethan Devinatz  has 
informed the author that all the results of \cite{DH} and \cite{LHS} go through as is, with $E_n$ 
replaced with $E(k, \Gamma)$ and $G_n$ replaced with 
$G(k) = S_n \rtimes \mathrm{Gal}(k/\mathbb{F}_p),$
a compact $p$-adic analytic group. Thus, \cite{DH} implies that 
(a) there is an isomorphism \[\pi_\ast(\lhat(E(k, \Gamma)^{dhK} \wedge E(k, \Gamma))) \cong \mathrm{Map}_c(G(k),\pi_\ast(E(k, \Gamma)))^K;\] and (b) $E(k, \Gamma)^{dhK} \simeq \lhat(\colim_j
E(k, \Gamma)^{dhU_jK}),$ where $\{U_j\}$ is a nested sequence of open normal subgroups of 
$G(k)$, as before.
\par
By (b), there is a continuous $G(k)$-spectrum 
\[E(k, \Gamma) \simeq \holim_i ((\colim_j E(k, \Gamma)^{dhU_j}) \wedge M_{i}),\] 
so that, for any closed subgroup $K$ of $G(k)$, $E(k, \Gamma)^{hK}$ exists. Then 
\[E(k, \Gamma)^{hK} \simeq E(k, \Gamma)^{dhK},\] by using the 
arguments of \cite[Theorem 8.2.1]{joint} and \cite{thesis}. Also, (a) and (b) 
imply that  $(\colim_j E(k, \Gamma)^{dhU_j}) \wedge M_{i}$ is a totally hyperfibrant 
discrete $G(k)$-spectrum, so that 
\[(E(k, \Gamma)^{hH})^{hK/H} \simeq E(k, \Gamma)^{hK},\] for any $H$ closed in $G(k)$ and normal 
in $K$.
\appendix
\section*{Appendix: An example of a discrete $G$-spectrum that is not hyperfibrant}
\begin{center}
\small\textsc{Daniel G. Davis, Ben Wieland}
\end{center}
\vspace{.1in}
\par
Let $G$ be a profinite group. In this appendix, 
we give an example, due to the second author, of a discrete $G$-spectrum $X$ that 
is not  a hyperfibrant discrete $G$-spectrum (for the meaning of this term, 
see Definition \ref{hyper}). All 
spectra are Bousfield-Friedlander spectra of simplicial sets.
\par
Let $p$ and $q$ be distinct primes, let 
$\mathbb{Z}_q = \lim_{n \geq 0} \mathbb{Z}/{q^n}$, and let 
\[G = \mathbb{Z}/p \times \mathbb{Z}_q.\] 
\par
For each $n \geq 0$, let $\mathbb{Z}/{q^n}$ act on itself. 
Denote by $\gpring$ the free $(\mathbb{Z}/{{p}})$-module on this 
set and by $H(\mathbb{Z}/p[\mathbb{Z}/{q^n}])$ the 
Eilenberg-MacLane spectrum of that abelian group. 
By the functoriality of these constructions, this is a 
$(\mathbb{Z}/{q^n})$-spectrum, for each $n$. 
By letting $\mathbb{Z}/{p}$ 
act trivially on $H(\mathbb{Z}/{{p}}[\mathbb{Z}/{q^n}])$, 
$H(\mathbb{Z}/{{p}}[\mathbb{Z}/{q^n}])$ is a (discrete) 
$(\mathbb{Z}/{{p}} \times \mathbb{Z}/{q^n})$-spectrum. 
Let $G$ act through the canonical surjection 
$G \rightarrow \mathbb{Z}/{{p}} \times \mathbb{Z}/{q^n}$, 
making $H(\mathbb{Z}/{{p}}[\mathbb{Z}/{q^n}])$ a discrete 
$G$-spectrum. 
By letting $G$ act trivially on $S^n$, the spectrum 
$\Sigma^n H(\mathbb{Z}/{{p}}[\mathbb{Z}/{q^n}])$ is also a discrete 
$G$-spectrum. Let \[H_n = \Hn\] and set 
\[X = \textstyle\bigvee_{n \geq 0} H_n.\] 
Since the colimit in the category of discrete $G$-spectra is formed in the category 
of spectra, $X$ is a discrete $G$-spectrum. 
We will show that $X$ is not a hyperfibrant discrete $G$-spectrum. In particular, 
we will show that the map $\psi(X)_{\mathbb{Z}/p}^G$ of (\ref{sigh}) is not a weak 
equivalence.
\par
For each $n$, $\pi_\ast(H_n)$ is 
$\mathbb{Z}/{{p}}[\mathbb{Z}/{q^n}]$ in degree $n$ and zero elsewhere, so that the 
$G$-equivariant map 
\[\phi \colon X = \textstyle\bigvee_{n \geq 0} H_n
\rightarrow \textstyle\prod_{n \geq 0} H_n,\] with 
$G$ acting diagonally on the target, 
is a weak equivalence. Thus, the induced map
\[X^{h\mathbb{Z}/p} \overset{\simeq}{\rightarrow} 
(\textstyle\prod_{n \geq 0} H_n)^{h\mathbb{Z}/p}\]
is a weak equivalence. 
\par
Let $H_n \rightarrow 
(H_n)_{f,\mathbb{Z}/p}$ be a 
trivial cofibration to a fibrant $(\mathbb{Z}/p)$-spectrum, all in the model category of 
(discrete) $(\mathbb{Z}/p)$-spectra. Then 
$(H_n)_{f,\mathbb{Z}/p}$ and $\textstyle\prod_{n \geq 0} (H_n)_{f,\mathbb{Z}/p}$ are 
fibrant spectra, and, hence, 
\[
(\underset{n \geq 0}{\textstyle\prod}H_n)^{h\mathbb{Z}/p} 
\simeq 
\displaystyle\holim_{\mathbb{Z}/p} \underset{n \geq 0}{\textstyle\prod} 
(H_n)_{f,\mathbb{Z}/p} 
\cong \underset{n \geq 0}{\textstyle\prod} \displaystyle\holim_{\mathbb{Z}/p} 
(H_n)_{f,\mathbb{Z}/p}
\simeq \underset{n \geq 0}{\textstyle\prod} 
H_n^{h\mathbb{Z}/p},\]
where the first equivalence follows, for example, from \cite[pg. 337]{cts}.
\par
It is well-known that $\cdp(\mathbb{Z}/p) = \infty.$ But, since $\cdq(\mathbb{Z}_q) = 1$ and 
$\mathbb{Z}_q$ is open in $G$, it follows that $G$ has finite virtual cohomological 
dimension. Thus, if $K$ is a closed subgroup of $G$, then, 
by \cite[Theorem 5.2]{fibrantmodel}, there is 
an identification 
\[X^{hK} = \holim_\Delta(\Gamma^\bullet_G(X_{f,G}))^K.\] There is a canonical map 
\[\psi \: \colim_{N \vartriangleleft_o G} X^{hN\mathbb{Z}/p} = 
\colim_{N \vartriangleleft_o G} \holim_\Delta(\Gamma^\bullet_G(X_{f,G}))^{N\mathbb{Z}/p} 
\rightarrow \holim_\Delta(\Gamma^\bullet_G(X_{f,G}))^{\mathbb{Z}/p} = X^{h\mathbb{Z}/p}.\] 
\par
Since $\mathbb{Z}/p \vartriangleleft G$, $N\mathbb{Z}/p \vartriangleleft G$, for each $N$. Hence, 
the above identification implies that $X^{hN\mathbb{Z}/p}$ is a 
$(G/{(N\mathbb{Z}/p}))$-spectrum, 
so that $X^{hN\mathbb{Z}/p}$ is also a discrete $G$-spectrum, since 
$G/{(N\mathbb{Z}/p)}$ is finite. This implies that $X^{hN\mathbb{Z}/p}$ is a discrete 
$\mathbb{Z}_q$-spectrum, so that $ \colim_{N \vartriangleleft_o G} X^{hN\mathbb{Z}/p}$, the 
source of $\psi$, is also a discrete $\mathbb{Z}_q$-spectrum. Notice that $X^{h\mathbb{Z}/p}$ 
is a $(G/(\mathbb{Z}/p))$-spectrum; that is, the target of $\psi$ is a $\mathbb{Z}_q$-spectrum. 
Also, it is easy to see that the map $\psi$ is $\mathbb{Z}_q$-equivariant. 
\par
Now suppose that $X$ is a hyperfibrant discrete $G$-spectrum. Thus, 
$\psi$ is a weak equivalence, so that $\pi_\ast(\psi)$ is a $\mathbb{Z}_q$-equivariant 
isomorphism between $\mathbb{Z}_q$-modules, such that the source is a \emph{discrete} 
$\mathbb{Z}_q$-module. This implies that $\pi_\ast(X^{h\mathbb{Z}/p})$ must also be a discrete 
$\mathbb{Z}_q$-module. Hence, to show that $X$ is not a hyperfibrant discrete $G$-spectrum, 
it suffices to show that $\pi_\ast(X^{h\mathbb{Z}/p})$ is not a discrete $\mathbb{Z}_q$-module. 
In particular, we only have to show that 
\[\pi_0(X^{h\mathbb{Z}/p}) \cong \textstyle\prod_{n \geq 0} 
\pi_0(H_n^{h\mathbb{Z}/p})\] fails to be a discrete 
$\mathbb{Z}_q$-module.  
\par
For each $n \geq 0$, since there is a descent spectral sequence 
\[E_2^{s,t} = H^s(\Zp; \pi_t(H_n)) \Rightarrow \pi_{t-s}(H_n^{h\mathbb{Z}/p}),\] 
where, for all $s$ and $t$,  
\[E_2^{s,t} =
\begin{cases}
H^s(\Zp; \gpring), & \mathrm{if} \ t=n; \\
0, & \mathrm{if} \ t \neq n,\end{cases}
\]  
\[\pi_t(H_n^{h\mathbb{Z}/p}) \cong 
\begin{cases}
H^{n-t}(\Zp; \gpring), & \mathrm{if} \, t \leq n; \\
0, & \mathrm{if} \ t > n.\end{cases}
\] In particular, 
\[\pi_0(H_n^{h\mathbb{Z}/p}) \cong 
H^n({\mathbb{Z}/{p}}\,; \mathbb{Z}/{{p}}[\mathbb{Z}/{q^n}]) \cong 
\bigoplus_{\mathbb{Z}/{q^n}} H^n(\mathbb{Z}/p \, ; \mathbb{Z}/p),\] 
since $\mathbb{Z}/{{p}}[\mathbb{Z}/{q^n}]$ is a free $(\Zp)$-module, with each 
copy of $\Zp$ having the trivial $(\mathbb{Z}/p)$-action. By applying 
\cite[Corollary 10.36]{Rotman}, it is not 
hard to see that \[H^n(\mathbb{Z}/p \, ; \mathbb{Z}/p) \cong \mathbb{Z}/p,\] for all $n \geq 0$, so 
that $\pi_0(H_n^{h\mathbb{Z}/p}) \cong 
\mathbb{Z}/{{p}}[\mathbb{Z}/{q^n}],$ for all $n \geq 0$. Hence, 
\[\pi_0(X^{h\mathbb{Z}/p}) \cong \textstyle\prod_{n \geq 0} \mathbb{Z}/{{p}}[\mathbb{Z}/{q^n}].\]  
\par
We now work to understand the $\Zq$-action on $X^{h\mathbb{Z}/p}$ so that we can show 
that $\pi_0(X^{h\mathbb{Z}/p})$ fails to be a discrete $\Zq$-module. 
If $K$ is any group and $Z$ is a $K$-spectrum, let $Z \rightarrow Z_\mathtt{f}$ be a weak 
equivalence to a fibrant spectrum; by 
the functoriality of fibrant replacement, $Z_\mathtt{f}$ is a $K$-spectrum and 
the weak equivalence is $K$-equivariant. If $Z$ is any spectrum, we 
let $\mathrm{Sets}(K, Z)$ be the spectrum that is 
defined by $\mathrm{Sets}(K,Z)_l = \mathrm{Sets}(K,Z_l),$ for each $l \geq 0$, where 
$\mathrm{Sets}(K,Z_l)$ is the pointed simplicial set with $k$-simplices equal to 
$\mathrm{Sets}(K,Z_{l,k})$, for all $k \geq 0$. Also, let $EK_\bullet$ be the usual free 
contractible simplicial $K$-set, with $EK_k = K^{k+1}$. 
\par 
Let 
$H_n^\mathtt{f} = (H_n)_\mathtt{f}$. Then, by \cite[proof of Lemma 10.5]{wilkerson}, 
\[H_n^{h\mathbb{Z}/p} 
= \holim_{[k] \in \Delta} \mathrm{Sets}(G^{k+1}, H_n^\mathtt{f})^{\Zp},\]
where the $(\Zp)$-fixed points are taken with respect to the $G$-action on 
the spectrum $\mathrm{Sets}(G^{k+1}, H_n^\mathtt{f})$ that is given by conjugation (that is, 
given $g \in G$ and $f$ in $\mathrm{Sets}(G^{k+1}, H_n^\mathtt{f})$, 
$(g \cdot f)((g_i)_i) = g \cdot f((g^{-1}g_i)_i),$ where $(g_i)_i \in G^{k+1}$ and, here, and elsewhere, 
we write maps in terms of elements because of the simplicial \emph{sets} that constitute all of our 
spectra). Since 
$\mathrm{Sets}(EG_\bullet, H_n^\mathtt{f})$ is a cosimplicial $G$-spectrum, 
$\holim_{[k] \in \Delta} \mathrm{Sets}(G^{k+1}, H_n^\mathtt{f})^{\Zp}$ is 
a $\Zq$-spectrum, showing that $H_n^{h\mathbb{Z}/p}$ is a $\Zq$-spectrum. 
\par
The $(\Zp)$-equivariant projection $G \rightarrow \Zp$ 
induces the weak equivalence
\[\lambda \colon \holim_{[k] \in \Delta} \mathrm{Sets}((\Zp)^{k+1}, H_n^\mathtt{f})^{\Zp} 
\rightarrow \holim_{[k] \in \Delta} \mathrm{Sets}(G^{k+1}, H_n^\mathtt{f})^{\Zp},\] 
where 
both the source and target of $\lambda$ are models for the $\Zq$-spectrum 
$H_n^{h\mathbb{Z}/p}$, with the $\Zq$-action on the source given by 
\[(g \cdot f)((h_i))_i = g \cdot f((h_i)_i),\] for $g \in \Zq$, 
$f \in \mathrm{Sets}((\Zp)^{k+1}, H_n^\mathtt{f})^{\Zp}$, and $(h_i)_i \in (\Zp)^{k+1}$. 
Given $h \in \Zp$,
\begin{align*} (h \cdot (g \cdot f))((h_i)) & = h \cdot ((g \cdot f)((h^{-1}h_i))) 
= hg \cdot f(h^{-1} \cdot (h_i)) 
\\ & = gh \cdot (h^{-1} \cdot f((h_i))) = (g \cdot f)((h_i)).
\end{align*}
Hence, $h \cdot (g \cdot f) = g \cdot f$, verifying that $g \cdot f \in 
\mathrm{Sets}((\Zp)^{k+1}, H_n^\mathtt{f})^{\Zp}$, as required.
\par
Now we show that $\lambda$ is $\Zq$-equivariant. Let $(h,g) \in G$, as above, 
and let 
\[\lambda_k \colon \mathrm{Sets}((\Zp)^{k+1}, H_n^\mathtt{f})^{\Zp} 
\rightarrow \mathrm{Sets}(G^{k+1}, H_n^\mathtt{f})^{\Zp}\] be the 
map, the collection of which induces $\lambda$. Explicitly, given $(h_i,g_i)_i \in G^{k+1}$,
\[\lambda_k(f)((h_i,g_i)) = f((h_i)).\] Notice that 
$\lambda_k(g \cdot f)((h_i,g_i)) = (g \cdot f)((h_i)) = g \cdot f((h_i))$ 
and 
\[(g \cdot \lambda_k(f))((h_i,g_i)) = g \cdot (\lambda_k(f)(g^{-1} \cdot (h_i, g_i))) 
= g \cdot (\lambda_k(f)((h_i,g^{-1}g_i))) = g \cdot f((h_i)),
\]
showing that $\lambda_k$ is $\Zq$-equivariant, and, hence, $\lambda$ is $\Zq$-equivariant.
\par 
To summarize, we have shown that 
the weak equivalence $\lambda$, 
between two different models for $H_n^{h\mathbb{Z}/p}$, 
is $\Zq$-equivariant. As will be seen, it will be convenient to use the source of 
$\lambda$, $\holim_{[k] \in \Delta} \mathrm{Sets}((\Zp)^{k+1}, H_n^\mathtt{f})^{\Zp},$ 
as our model for the $\Zq$-spectrum $H_n^{h\mathbb{Z}/p}$.
\par
Since $\phi$ induces the composition 
$(\textstyle\bigvee_{n \geq 0} H_n)_\mathtt{f} \rightarrow (\prod_{n \geq 0} H_n)_\mathtt{f} 
 \rightarrow (\prod_{n \geq 0} H_n^\mathtt{f})_\mathtt{f}$, which is 
$G$-equivariant and a weak equivalence, there are weak equivalences 
\begin{align*}
X^{h\mathbb{Z}/p} & = 
\holim_{[k] \in \Delta} \mathrm{Sets}(G^{k+1}, 
(\underset{n \geq 0}{\textstyle\bigvee} H_n)_\mathtt{f})^{\Zp} 
\overset{\simeq}{\rightarrow} \displaystyle\holim_{[k] \in \Delta} 
\mathrm{Sets}(G^{k+1}, 
(\underset{n \geq 0}{\textstyle\prod} H_n^\mathtt{f})_\mathtt{f})^{\Zp} \\ 
& \overset{\simeq}{\leftarrow} \displaystyle\holim_{[k] \in \Delta} 
\mathrm{Sets}(G^{k+1}, 
\underset{n \geq 0}{\textstyle\prod} H_n^\mathtt{f})^{\Zp} 
\cong \underset{n \geq 0}{\textstyle\prod} \displaystyle\holim_{[k] \in \Delta} 
\mathrm{Sets}(G^{k+1}, H_n^\mathtt{f})^{\Zp} \\
& \overset{\simeq}{\leftarrow} \underset{n \geq 0}{\textstyle\prod} 
\displaystyle\holim_{[k] \in \Delta} \mathrm{Sets}((\Zp)^{k+1}, 
H_n^\mathtt{f})^{\Zp} = \underset{n \geq 0}{\textstyle\prod} H_n^{h\mathbb{Z}/p},
\end{align*}
where the last weak equivalence uses that $\lambda$ is a weak equivalence. All of the maps in the above zigzag are $\Zq$-equivariant, so that 
\[\pi_0(X^{h\mathbb{Z}/p}) \cong \textstyle\prod_{n \geq 0} \pi_0(H_n^{h\mathbb{Z}/p})\] 
is a $\Zq$-equivariant isomorphism, with $\Zq$ acting diagonally on 
$\textstyle\prod_{n \geq 0} \pi_0(H_n^{h\mathbb{Z}/p})$.
\par
From the descent spectral sequence for $H_n^{h\mathbb{Z}/p}$, we have
\begin{align*}
\pi_0(H_n^{h\mathbb{Z}/p}) & \cong E_2^{n,n} 
= H^{n}[\,\mathrm{Sets}((\Zp)^{\ast+1}, \pi_n(H_n))^{\Zp}\,] \\
& \cong H^n[\,\mathrm{Sets}((\Zp)^{\ast+1}, \gpring)^{\Zp}\,],
\end{align*} where $\Zp$ acts trivially on $\gpring$ and $\Zq$ acts on 
$\gpring$ through the canonical map $\Zq \rightarrow \mathbb{Z}/q^n$ and the 
usual action of $\mathbb{Z}/q^n$ on $\gpring$. (We are saying that 
$\mathbb{Z}/q^n$ acts on $\pi_n(\Sigma^n H(\gpring))$ in exactly the same way that 
$\mathbb{Z}/q^n$ acts on $\gpring$, at the beginning of our argument. 
That this is true follows from the proof of \cite[Proposition 3.4]{hGal} and 
the discussion in \cite{hGal} between Lemma 3.1 and Proposition 3.4.) Thus, given 
$g \in \Zq$ and $(r_j j)_j \in \gpring$, where $r_j \in \Zp$ and 
$j$ runs through the elements of $\mathbb{Z}/q^n$, $g \cdot (r_j j)_j = (r_j gj)_j.$
\par
Notice that in
\begin{align*}
\pi_0(H_n^{h\mathbb{Z}/p}) & \cong
H^n[\,\mathrm{Sets}((\Zp)^{\ast+1}, \gpring)^{\mathbb{Z}/p}\,] \\
& \cong H^n[\,\textstyle\prod_{\mathbb{Z}/q^n} 
\mathrm{Sets}((\Zp)^{\ast+1}, \Zp)^{\mathbb{Z}/p}\,],
\end{align*} the first isomorphism is 
$\Zq$-equivariant. Since we require the second isomorphism to also be 
$\Zq$-equivariant, we now seek to understand why this is the case.
\par
If $f$ is in $\mathrm{Sets}((\Zp)^{k+1}, 
\gpring)^{\mathbb{Z}/p},$ then 
\[f(h) = (f_j(h) j)_j,\] where 
$f_j \in \mathrm{Sets}((\Zp)^{k+1}, \Zp)$. If $\eta \in \Zp$, then 
\[(f_j(h)j)_j = f(h) = (\eta \cdot f)(h) = 
\eta \cdot (f_j(\eta^{-1}h)j)_j = (f_j(\eta^{-1}h)j)_j,\] so that 
$f_j(\eta^{-1}h) = f_j(h) = \eta^{-1} \cdot f_j(h),$ since $\Zp$ acts trivially on 
$f_j(h) \in \Zp$. Thus, $f_j(h) = \eta \cdot f_j(\eta^{-1}h) = (\eta \cdot f_j)(h),$
so that $f_j \in \mathrm{Sets}((\Zp)^{k+1}, \Zp)^{\mathbb{Z}/p}.$ 
\par
The preceding conclusion implies that the isomorphism
\[\theta \colon \mathrm{Sets}((\Zp)^{k+1}, 
\gpring)^{\mathbb{Z}/p} \rightarrow 
\textstyle\prod_{\mathbb{Z}/q^n} 
\mathrm{Sets}((\Zp)^{k+1}, \Zp)^{\mathbb{Z}/p}\] is given by 
\[\theta(f) = (f_jj)_j,\] so that the $j$th coordinate of $\theta(f)$ is $f_j$. 
Notice that, given $g \in \Zq$, 
\[(g \cdot f)(h) = (f_j(h) gj)_j = (f_{g^{-1}j}(h)j)_j.\] 
Thus, $\theta(g \cdot f) = (f_{g^{-1}j} j)_j.$ 
\par
Since we require $\theta$ to be $\Zq$-equivariant, we must have 
\[g \cdot \theta(f) = g \cdot (f_j j)_j = (f_{g^{-1}j}j)_j.\] 
We conclude that if $(k_jj)_j \in \prod_{\mathbb{Z}/q^n} 
\mathrm{Sets}((\Zp)^{k+1}, \Zp)^{\mathbb{Z}/p}$, then the $\Zq$-action is 
defined by 
\[g \cdot (k_jj)_j = (k_{g^{-1}j}j)_j;\] it is easy to verify that this is indeed 
an action. Therefore, we have shown that there is 
a $\Zq$-equivariant isomorphism of $\Zq$-modules
\[\pi_0((H_n)^{h\mathbb{Z}/p}) \cong H^n[\,\textstyle\prod_{\mathbb{Z}/q^n} 
\mathrm{Sets}((\Zp)^{\ast+1}, \Zp)^{\mathbb{Z}/p}\,],\] 
where the $\Zq$-action on the right-hand side is induced by 
$g \cdot (k_jj)_j = (k_{g^{-1}j}j)_j,$ given 
$(k_jj)_j \in \prod_{\mathbb{Z}/q^n} 
\mathrm{Sets}((\Zp)^{k+1}, \Zp)^{\mathbb{Z}/p}$. It is useful to note the following 
feature of this $\Zq$-action: when $g$ acts on $(k_jj)_j$ to yield $(k_{g^{-1}j}j)_j$, 
the action is only moving around the coordinate functions $k_j$, but it is not 
changing them. 
\par
The aforementioned ``feature" implies that, in the isomorphisms 
\begin{align*}
H^n[\,\textstyle\prod_{\mathbb{Z}/q^n} 
\mathrm{Sets}((\Zp)^{\ast+1}, \Zp)^{\mathbb{Z}/p}\,]  
& \cong 
\textstyle\prod_{\mathbb{Z}/q^n} 
H^n[\,\mathrm{Sets}((\Zp)^{\ast+1}, \Zp)^{\mathbb{Z}/p}\,] \\
& \cong \textstyle\prod_{\mathbb{Z}/q^n} H^n(\Zp; \Zp) \cong \textstyle\prod_{\mathbb{Z}/q^n} \Zp \\
& \cong \gpring,
\end{align*}
the $\Zq$-action on the first term  
$H^n[\,\textstyle\prod_{\mathbb{Z}/q^n} 
\mathrm{Sets}((\Zp)^{\ast+1}, \Zp)^{\mathbb{Z}/p}\,]$ corresponds to the 
following action on the last term $\gpring$: 
\[g \cdot (r_j j)_j = (r_{g^{-1}j}j)_j = (r_j gj)_j.\] Therefore, we have shown that 
$\pi_0(X^{h\mathbb{Z}/p}) \cong \textstyle\prod_{n \geq 0} \gpring$
is a $\Zq$-equivariant isomorphism of $\Zq$-modules, where the $\Zq$-action on the 
right-hand side is defined as follows: given $g \in \Zq$ and an element 
$((r_{j_n} j_n)_{j_n})_{n \geq 0}$ in $\textstyle\prod_{n \geq 0} \gpring,$ 
where $(r_{j_n}j_n)_{j_n} \in \gpring,$ 
\[g \cdot ((r_{j_n}j_n)_{j_n})_{n \geq 0} = ((r_{j_n}gj_n)_{j_n})_{n \geq 0}.\] It only 
remains to show that $\prod_{n \geq 0} \gpring$ is not a discrete $\Zq$-module.
\par
Consider the injection of sets
\[i \colon \Zq = \lim_{n \geq 0} \mathbb{Z}/q^n \rightarrow 
\textstyle\prod_{n \geq 0} \gpring, 
\ \ \ (\alpha_n)_{n \geq 0} \mapsto ((i(\alpha_n,j_n)j_n)_{j_n})_{n \geq 0},\]
where $\alpha_n \in \mathbb{Z}/q^n$ and 
\[i(\alpha_n, j_n) = \begin{cases}
1, & \mathrm{if} \ j_n = \alpha_n; \\
0, & \mathrm{if} \ j_n \neq \alpha_n.
\end{cases}
\] 
\par
Note that $i$ is $\Zq$-equivariant. We have that 
\[i(g \cdot (\alpha_n)) = i((g\alpha_n)) = ((i(g\alpha_n,j_n)j_n))\] and 
\[g \cdot i((\alpha_n)) = g \cdot ((i(\alpha_n,j_n)j_n)) = ((i(\alpha_n,j_n)gj_n)_{j_n}).\] 
The only nonzero coordinate of $(i(g\alpha_n,j_n)j_n) \in \gpring$ 
is the $(g\alpha_n)$th 
coordinate, which is equal to $1 g\alpha_n$. Similarly, the only nonzero 
coordinate of the tuple $(i(\alpha_n,j_n)gj_n)_{j_n}$ is $1gj_n$, where $j_n=\alpha_n$, 
yielding the element $1g\alpha_n$, the $(g\alpha_n)$th coordinate. Thus, 
$i(g \cdot (\alpha_n)) = g \cdot i((\alpha_n))$, so that $i$ is indeed $\Zq$-equivariant. 
\par
Since the function $i$ is a 
$\Zq$-equivariant injection, we can regard the $\Zq$-set 
$\Zq$ as a $\Zq$-subset of $\prod_{n \geq 0} \gpring$. Now suppose that 
$\prod_{n \geq 0} \gpring$ is a discrete $\Zq$-module. Let $\GG$ be the profinite 
group $\Zq$ and let the $\GG$-set $\Zq$ have the discrete topology. Then the composition 
\[\GG \times \Zq \overset{\mathrm{id} \times i}{\longrightarrow} 
\GG \times (\textstyle\prod_{n \geq 0} \gpring) 
\overset{\mathrm{action}}{\longrightarrow} \textstyle\prod_{n \geq 0} \gpring\] is continuous, and, since the 
image of this map is inside $\Zq$, the induced 
action map $\GG \times \Zq \rightarrow \Zq$
is continuous, so that $\Zq$ is a discrete $\GG$-set. But this is 
a contradiction, since $\Zq$ is a profinite $\GG$-space and not a discrete $\GG$-set. 
Therefore, $\prod_{n \geq 0} \gpring$ is not a discrete $\Zq$-module. 
\par
We have completed the proof that $X$ is not a hyperfibrant discrete $G$-spectrum. 
The same construction also works in the case that $p=q$, but the 
analysis is a little more complicated. This realizes a suggestion of the 
referee of this paper, who had suggested  to the first 
author, independently of the second author's work, 
that there might exist a non-hyperfibrant discrete $(\mathbb{Z}/p \times \mathbb{Z}_p)$-spectrum.

\bibliographystyle{plain}
\bibliography{biblio}
\vspace{.25in}

\end{document}